\documentclass[12pt]{amsart}
\usepackage{amssymb,amscd}
\usepackage[all,ps,cmtip]{xy}
\makeatletter
\def\cl@part{}
\@addtoreset{equation}{part}
\@addtoreset{equation}{section}
\def\swappedhead#1#2#3{%
(\thmnumber{#2})\thmname{\@ifnotempty{#2}{\ }#1}%
\thmnote{ {\the\thm@notefont(#3)}}}
\makeatother

\newtheorem{theoreme}[equation]{Theorem}
\newtheorem{proposition}[equation]{Proposition}
\newtheorem{lemme}[equation]{Lemma}
\newtheorem{corollaire}[equation]{Corollary}

\theoremstyle{definition}
\newtheorem{definition}[equation]{Definition}

\theoremstyle{remark}
\newtheorem{rem}[equation]{Remark}
\newtheorem{rems}[equation]{Remarks}
\newtheorem*{rem*}{Remark}
\newtheorem*{rems*}{Remarks}

\DeclareMathOperator{\GL}{GL}

\DeclareMathOperator{\spe}{sp}
\DeclareMathOperator{\spec}{Spec}
\DeclareMathOperator{\supp}{Supp}
\DeclareMathOperator{\Hom}{Hom}
\DeclareMathOperator{\Ext}{Ext}

\DeclareMathOperator{\gr}{gr}
\DeclareMathOperator{\car}{Char}
\DeclareMathOperator{\SP}{SP}

\DeclareMathOperator{\can}{can}

\DeclareMathOperator{\DR}{DR}
\renewcommand{\ker}{\mathop{\rm Ker}\nolimits}
\DeclareMathOperator{\rk}{rk}
\let\rg\rk
\DeclareMathOperator{\tr}{tr}
\DeclareMathOperator{\ir}{ir}
\DeclareMathOperator{\id}{Id}
\DeclareMathOperator{\im}{Im}
\DeclareMathOperator{\image}{image}
\def\dtau{\partial_\tau}
\def\dth{\partial_\theta}
\def\pbR{{{}^p\!\bR}}
\def\pQQ{{{}^p\!\QQ}}
\def\pCC{{{}^p\!\CC}}
\def\psip{{{}^p\!\psi}}
\def\phip{{{}^p\!\phi}}
\def\Ct{\CC[\tau,\tau_{}^{-1}]}
\def\ct{\CC[\tau,\tau_{}^{-1}]}

\def\pDR{{{}^p\!\DR}}
\def\DRa{{\DR}^{\rm an}}
\def\pDRa{\pDR^{\rm an}}

\def\Afua{\AA\!^{1\rm an}}
\def\afua{\AA\!^{1\rm an}}
\def\Afu{\AA\!^1}
\def\afu{{\AA}\!^1}
\def\Afuc{\check\AA\!^1}
\def\afuc{\check\AA\!^1}

\let\PP\PP
\def\Pu{\PP^1}

\def\pch{{\check p}}

\def\Tc{\check T}

\def\Etf{\cE^{-\tau f}}

\def\psim{\psi^{\rm mod}}
\def\phim{\phi^{\rm mod}}

\def\muM{\cM^{\!{}^\mu}}

\newcommand{\ooplus}{\mathop\oplus\limits}
\newcommand{\ootimes}{\mathop\otimes\limits}
\newcommand{\lefpar}{\left(}
\newcommand{\rigpar}{\right)}
\newcommand{\lefcro}{\left[}
\newcommand{\rigcro}{\right]}
\newcommand{\lefacc}{\left\{}

\newcommand{\rig}{\right.}
\newcommand{\lcr}{\left[\!\left[}
\newcommand{\rcr}{\right]\!\right]}

\newcommand{\norme}[1]{\left\Vert#1\right\Vert}
\let\dpl\displaystyle
\let\wh\widehat
\let\wt\widetilde
\let\ov\overline
\let\wwwh\wh
\let\wwh\wh
\def\implique{\Longrightarrow}

\newcommand{\defin}{:=}
\newcommand{\apriori}{{\em a priori\/}}
\newcommand{\loccit}{{\em loc.\ cit}}
\newcommand{\eg}{{\it e.g}}
\newcommand{\cf}{\emph{cf}}

\newcommand{\ie}{{\it i.e}}
\newcommand{\resp}{{\it resp}}
\newcommand{\T}{\S\kern .15em }
\newcommand{\ptbl}{.\kern .15em }
\newcommand{\bbullet}{{\scriptscriptstyle\bullet}}

\DeclareMathAlphabet{\mathcalmaigre}{U}{eus}{m}{n}

\def\AA{\mathbb{A}}
\def\CC{\mathbb{C}}
\def\DD{\mathbb{D}}
\def\ZZ{\mathbb{Z}}
\def\QQ{\mathbb{Q}}
\def\NN{\mathbb{N}}

\def\PP{\mathbb{P}}

\def\cC{\mathcal{C}}
\def\cD{\mathcal{D}}
\def\cE{\mathcal{E}}
\def\cF{\mathcal{F}}
\def\cG{\mathcal{G}}
\def\cH{\mathcal{H}}

\def\cK{\mathcal{K}}
\def\cL{\mathcal{L}}
\def\cM{\mathcal{M}}
\def\cN{\mathcal{N}}
\def\cO{\mathcal{O}}

\def\cR{\mathcal{R}}
\def\cS{\mathcal{S}}

\def\cU{\mathcal{U}}

\def\cX{\mathcalmaigre{X}}
\def\cY{\mathcal{Y}}
\def\cZ{\mathcal{Z}}
\def\ccF{\mathcalmaigre{F}}
\def\ccG{\mathcalmaigre{G}}
\def\ccT{\mathcalmaigre{T}}
\def\ccV{\mathcalmaigre{V}}

\let\gG\gotG

\let\gM\gotM

\def\gammag{\boldsymbol{\gamma}}

\def\varepsilong{\boldsymbol{\varepsilon}}

\def\bmg{\boldsymbol{g}}

\def\bD{\boldsymbol{D}}

\def\bH{\boldsymbol{H}}

\def\bL{\boldsymbol{L}}

\def\bR{\boldsymbol{R}}

\newdimen\lengtharrow
 \newbox\exponantbox \newbox\indicebox
\def\dimmax#1#2{\ifdim#1<#2 #2\else #1\fi}
\def\arrowr#1#2%
    {\setbox\exponantbox=\hbox{$#2$}
      \setbox\indicebox=\hbox{$#1$}
       \lengtharrow=\dimmax{\wd\indicebox}{\wd\exponantbox}
         \ifdim \lengtharrow=0pt	
           \mathrel{\smash{\mathop{\hbox to 10mm{\rightarrowfill}}}}
             \else \ifdim\lengtharrow<6mm
               \lengtharrow=10mm
                 \else \advance\lengtharrow by 4mm
                   \fi
\mathrel{\smash{\mathop{\hbox to\lengtharrow{\rightarrowfill}}\limits%
                         ^{\textstyle#2}_{\textstyle#1}}}
\fi}
\def\arrowl#1#2%
    {\setbox\exponantbox=\hbox{$#2$}
     \setbox\indicebox=\hbox{$#1$}
       \lengtharrow=\dimmax{\wd\indicebox}{\wd\exponantbox}
         \ifdim \lengtharrow=0pt
           \mathrel{\smash{\mathop{\hbox to 10mm{\leftarrowfill}}}}
             \else \ifdim\lengtharrow<6mm
               \lengtharrow=10mm
                 \else \advance\lengtharrow by 4mm
                   \fi
            \mathrel{\smash{\mathop{\hbox to\lengtharrow{\leftarrowfill}}\limits%
                         ^{\textstyle#2}_{\textstyle#1}}}
                         \fi}
\def\hookarrowr#1#2{\lhook\joinrel\arrowr{#1}{#2}}

\def\arrowd#1#2{\llap{$ \textstyle #1$}\left\downarrow%
\vbox to 6mm{}\right.\rlap{$\textstyle #2$}}
\def\arrowu#1#2{\llap{$\textstyle #1$}\left\uparrow%
\vbox to 6mm{}\right.\rlap{$\textstyle #2$}}

\def\MRE#1{\arrowr{}{#1}}

\def\HMRE#1{\hookarrowr{}{#1}}

\def\MDR#1{\arrowd{}{#1}}
\def\MDL#1{\arrowd{#1}{}}
\def\MUR#1{\arrowu{}{#1}}

\newcommand{\isom}{\stackrel{\sim}{\longrightarrow}}

\begin{document}

\title{Hypergeometric periods for a tame polynomial}
\author{Claude Sabbah}
\address{URA 169 du C.N.R.S., Centre de Math\'ematiques, \'Ecole polytechnique, F--91128 Palaiseau cedex, France}
\email{sabbah@math.polytechnique.fr}
\begin{abstract}
We analyse the Gauss-Manin system of differential equations---and its Fourier transform---attached to regular functions satisfying a tameness assupmption on a smooth affine variety over $\mathbb{C}$ (e.g. tame polynomials on $\mathbb{C}^{n+1}$). We give a solution to the Birkhoff problem and prove Hodge-type results analogous to those existing for germs of isolated hypersurface singularities.
\end{abstract}
\date{June 1997}
\maketitle
\tableofcontents

\section*{Introduction}

Let $f:\CC_{}^{n+1}\rightarrow \CC$ be a nonconstant polynomial function. The hypergeometric periods associated with this polynomial are the integrals
\begin{eqnarray*}
I_{\gamma,\omega}^{}(s)&=&\int_\gamma f^s\omega
\end{eqnarray*}
where $\omega$ is an algebraic $(n+1)$-form and $\gamma$ an $(n+1)$-cycle with coefficients in a suitable local system. In \cite{Varchenko89} A.~N.~Varchenko gave a formula for the determinant of a period matrix made with such integrals for some specific polynomials: this determinant is expressed as a product of terms of the form $\Gamma(s+\beta)$ where $\beta$ varies in the {\em spectrum} of the polynomial (for these polynomials, there exists only one critical point and the spectrum is the spectrum of this singularity in the sense of \cite{Varchenko82}); moreover, in this situation, the choice of the forms $\omega$ is quite natural.

In \cite{Douai93}, A.~Douai considered a spectrum for convenient nondegenerate polynomials (in terms of the Newton filtration) and conjectured that there should exist a basis of forms $\omega$ so that the same formula holds (up to periodic functions in $s$) for the determinant of the period matrix. He proved this conjecture for some special cases.

The appearance of a product of $\Gamma$ factors is explained by the fact that such a determinant (whatever the choice of the forms $\omega$ or the cycles $\gamma$ can be) is a solution of a system of linear finite difference equations and one expects (if it is nonzero) that it is expressed (up to a periodic function in $s$) as a product of $\Gamma(s+\beta)$ where $\beta$ is a logarithm of an eigenvalue of the monodromy of $f$ at infinity \cite{L-S90} (we assume here that the monodromy around $0$ is trivial, see below). Changing the forms will change the $\beta$'s by an integer.

The problem addressed in this paper consists in finding a natural choice of such logarithms, called the spectrum, and in showing that there exists a family of differential forms such that the previously aluded results still hold. We will assume that the polynomial is tame (cohomologically tame will be enough, see \T\ref{secproptame}), \ie. that, for some compactification of $f$, no modification of the topology (or the cohomology) of the fibres comes from infinity.

For a cohomologically tame polynomial, there are two possible ways of defining a spectrum.

The first one uses the Jacobian quotient $\CC[x_0,\ldots ,x_n]/(\partial_{x_0}^{}f,\ldots ,\partial_{x_n}^{}f)$ as Var\-chenko does for isolated hypersurfaces singularities (refered to as the local case below). One should moreover use a filtration which measures an asymptotic behaviour, as in the local case. In the convenient nondegenerate case, Douai uses the Newton filtration. It turns out that, in general, one should use the filtration measuring the asymptotic behaviour of integrals $\int_\delta\omega e_{}^{-\tau f}$ where $\delta$ is a Lefschetz thimble (see \cite{Pham81}) and $\omega$ as above, when $\tau\rightarrow 0$ (and not when $\tau\rightarrow \infty$ as is usually done in the stationary phase method).

This can be translated into more algebraic terms. Let $M$ be the Gauss-Manin system of the polynomial $f$: this is a (regular holonomic) module on the Weyl algebra $\CC[t]\langle\partial_t\rangle$, which associated local system (outside the critical values of $f$) is the one made by the spaces $H^n(f_{}^{-1}(t),\CC)$. The {\em Brieskorn lattice} $M_0$ is a free $\CC[t]$-module inside $M$ defined as in the local situation (see \eg. \cite{Pham79}). It turns out that $\partial_t$ acts in a one-to-one way on $M$ and that $M_0$ is stable under $\partial_{t}^{-1}$ (as in the local case). The role of microlocalization in the local case is now played by the Fourier-Laplace transform (see however \cite{Malgrange81} for a direct interpretation of the microlocalization in terms of Fourier transform in the local case).

Let $G$ be the module $M$ when viewed as a $\CC[\tau]\langle\partial_\tau\rangle$-module, with $\tau=\partial_t$ and $\partial_\tau=-t$. This is a holonomic module on the affine line $\Afuc$ with coordinate $\tau$, which has singularities at $\tau=0$ and $\tau=\infty$ only, the former being regular, but not the latter in general (see \eg. \cite{Malgrange91}). Consider the coordinate $\theta=\tau_{}^{-1}=\partial_{t}^{-1}$ at infinity on the affine line $\Afuc$. Then $M_0$ is also a free $\CC[\theta]$-module (with an action of $t$), and we denote it $G_0$. The fibre at $\theta=0$ of this module (\ie. $G_0/\theta G_0$) is identified with the Jacobian quotient, up to the choice of a volume form. We are interested in asymptotic expansions of sections of $G_0$ when $\tau\rightarrow 0$. More precisely we consider the Malgrange-Kashiwara filtration $V_\bbullet G$ at $\tau=0$ and the spectrum corresponds to the jump indices of this filtration, when induced on the Jacobian quotient.

Once such a definition is given (it coincides with the one given by Douai in the convenient nondegenerate case, see \T\ref{secnondeg}), it remains to show the existence of good differential forms. We argue in two steps.

The hypergeometric determinant is closely related to the determinant of the Aomoto complex \cite{L-S90}
\begin{eqnarray*}
0\rightarrow \CC(s)\ootimes_\CC \Omega^0[1/f]\MRE{d_s}\cdots \MRE{d_s}\CC(s)\ootimes_\CC\Omega_{}^{n+1}[1/f]\rightarrow 0
\end{eqnarray*}
where $\Omega^k[1/f]$ denotes the space of rational $k$-forms on $\CC_{}^{n+1}$ with poles along $f=0$ and $d_s$ is the twisted differential $d_s=f_{}^{-s}\cdot d\cdot f^s$. The $(n+1)$th cohomology group of this complex (which is also the only possible nonzero one, thanks to tameness) is equal to the Mellin transform of the Gauss-Manin system $M$. It is also equal to the Mellin transform of the Fourier transform of $M$. When $0$ is not a critical value for $f$, a good family of differential forms will be obtained from a {\em good basis} of the Brieskorn lattice $G_0$.

This notion of a good basis has been introduced in the local case by M. Saito \cite{MSaito89} in order to show the existence of primitive forms conjectured by K. Saito. We give the straightforward adaptation to our situation, so that the main question which remains is the existence of such a good basis. M. Saito also gave a criterion for the existence of such a basis: the filtration induced by $G_k=\tau^kG_0$ on the nearby cycles of $G$ at $\tau=0$ should be the Hodge filtration of a mixed Hodge structure, the weight filtration of which should be decreased by two under the action of the nilpotent part of the monodromy of $G$.

This leads to the second possible definition of the spectrum: J. Steenbrink and S. Zucker \cite{S-Z85} constructed a limit mixed Hodge structure on the cohomology $H^n(f_{}^{-1}(t),\CC)$ when $t\rightarrow \infty$ (this construction is also given by F. El~Zein \cite{Elzein86} and M. Saito in \cite{MSaito87}), and using the procedure defined by J. Steenbrink in \cite{Steenbrink87}, one defines a spectrum, called the {\em spectrum at infinity} of the polynomial $f$. However it has no \apriori\ direct link with the spectrum of the Brieskorn lattice.

The main result of \cite{Bibi96a} consists in showing the existence of a natural mixed Hodge structure on the vanishing cycles at $\tau=0$ of the Fourier transform $G$ isomorphic to the mixed Hodge structure above. In particular both give the same spectrum.

In part II of the present paper we show that the filtration given by the Brieskorn lattice coincides with the Hodge filtration constructed in \cite{Bibi96a}, when the polynomial is cohomologically tame. The ideas are very similar to the one of A.~N.~Varchenko in the local case (see also \cite{Pham83}, \cite{S-S85}, \cite{MSaito89}). In particular we show that the spectrum at infinity and the spectrum of the Brieskorn lattice coincide (up to a shift for the eigenvalue $1$ of the monodromy), and this gives the existence of good bases, thanks to the criterion of M. Saito.

It should be noticed, as a consequence of the existence of such bases, that the Riemann-Hilbert-Birkhoff problem (\ie. the existence of a Birkhoff normal form) can be solved for the Brieskorn lattice: it is possible to find a trivial bundle on the projective line associated with $\Afuc$, which is contained in $G$, which restricts to $G_0$ in the chart $\tau\neq0$ and which restricts to a logarithmic connection with pole at $\tau=0$ on $\Afuc$.

\medskip
Part I gives sufficient conditions (propositions \ref{propbbases} and \ref{propbbasesdual}, following M. Saito) to solve the Birkhoff problem for some germs of irregular meromorphic connections in one variable, and makes the link with Mellin transform.

Part II is concerned with applications (see \T\ref{secappli}) to the Gauss-Manin system of a cohomologically tame polynomial on $\CC_{}^{n+1}$: to prove the existence of a good basis for this meromorphic connection, we apply a criterion of M. Saito (see \T\ref{rembbases}); in order to do so, we need to establish Hodge properties for this system, and for this we use the results of \cite{Bibi96a}. In fact it should be emphasized that all the results of part II are given for a regular function $f:U\rightarrow \Afu$ on an affine manifold, this function satisfying a tameness assumption at infinity on $U$. The only difference with the case $U=\AA_{}^{n+1}$ is that $\partial_t$ is not necessarily bijective on $M$, so one has to distinguish between $M_0$ and $G_0$. In particular, we solve the Birkhoff problem for $G_0$ and find a good basis for $G_0$, but this is not directly translated as properties on $M_0$.

\medskip
We refer to the appendix of \cite{Bibi96a} for the notation which is not defined in the main course, in particular for the conventions made concerning perverse functors.

\medskip
We thank A. Douai, R. Garcia and F. Loeser for useful discussions.

\part{Good basis of a meromorphic connection \mbox{with a lattice}}

\section{Spectrum of a meromorphic connection with a lattice}\label{sec1}
We shall denote $U_0=\spec\CC[\tau]$ and $U_\infty=\spec\CC[\theta]$ the two standard charts of $\PP^1(\CC)$, where $\theta=1/\tau$ on $U_0\cap U_\infty$. We shall denote $0=\{\tau=0\}$ and $\infty =\{\theta=0\}$.

Let $G$ be a meromorphic connection on $\PP^1$, with singularities at $0$ and $\infty$ only, the singularity at $0$ being regular (but not necessarily the one at infinity). Then $G$ is a free $\Ct$-module of finite rank $\mu$, equipped with a derivation $\dtau$ which makes it a left $\Ct\langle\dtau\rangle$-module. In the following, we shall identify $\Ct$ with $\CC[\theta,\theta_{}^{-1}]$ by $\theta=1/\tau$.

\subsection*{Lattices}
A {\em lattice} of $G$ is a free $\CC[\theta]$-submodule $G_0$ of $G$ such that $\CC[\theta,\theta_{}^{-1}]\otimes_{\CC[\theta]}G_0=G$. Such a lattice $G_0$ has then rank $\mu$ over $\CC[\theta]$. We will say that the lattice $G_0$ has {\em type $1$ at infinity} if moreover $G_0$ is stable under the action of the operator $t\defin\theta^2\dth$.

In this paper, all lattices considered will have type $1$ at infinity, so we will simply call them lattices. We shall denote $G_k=\theta_{}^{-k}G_0$, in order to obtain an increasing filtration of $G$ by lattices.

\subsection*{Malgrange-Kashiwara filtration of $G$ at $0$}
Consider the connection $G$ on $U_0$. It is known that it is a holonomic $\CC[\tau]\langle\dtau\rangle$-module. Its only singular point is $0$. In order to simplify the argument, we shall assume in the following that the monodromy of $\DRa G$ on $U_0-\{0\}$ is quasi-unipotent. In fact, we shall only apply what follows to this situation.

We shall now use in this simple situation the properties of the Malgrange-Kashiwara filtration (see \eg. \cite{Bibi87,M-S86,MSaito86}, or \cite[\T6]{Bibi90} for the one dimensional case) that we briefly recall. Let $V_\bbullet\CC[\tau]\langle\dtau\rangle$ be the increasing filtration of $\CC[\tau]\langle\dtau\rangle$ defined by
\begin{eqnarray*}
V_{-k}\CC[\tau]\langle\dtau\rangle&=&\tau^k\CC[\tau]\langle\tau\dtau\rangle\quad\mbox{for }k\geq 0\\
V_{k}\CC[\tau]\langle\dtau\rangle&=&V_{k-1}\CC[\tau]\langle\tau\dtau\rangle+\dtau V_{k-1}\CC[\tau]\langle\tau\dtau\rangle\quad\mbox{for }k\geq 1.
\end{eqnarray*}

There exists a unique increasing exhaustive filtration $V_\bbullet G$ of $G$, indexed by the union of a finite number of subsets $\alpha+\ZZ\subset\QQ$, satisfying the following properties:
\begin{enumerate}
\item
For every $\alpha$, the filtration $V_{\alpha+\ZZ}G$ is good relatively to $V_\bbullet\CC[\tau]\langle\dtau\rangle$;
\item
For every $\beta\in\QQ$, $\tau\dtau+\beta$ is nilpotent on $\gr_{\beta}^{V}G\defin V_\beta G/V_{<\beta}G$.
\end{enumerate}

By assumption, each $V_\beta G$ is a finite type module over $\CC[\tau]\langle\tau\dtau\rangle$. Because $\tau$ is invertible on $G$ the map induced by $\tau$
\begin{eqnarray*}
\tau:V_\beta G\longrightarrow V_{\beta-1}G
\end{eqnarray*}
is bijective. Consequently we have for every $\beta\in\QQ$
\begin{eqnarray}\label{eqnresb}
\Ct\ootimes_{\CC[\tau]}V_\beta G=G.
\end{eqnarray}

Let $G_{}^{\rm an}=\CC\{\tau\}\otimes_{\CC[\tau]}G$ and $VG_{}^{\rm an}$ be the Malgrange-Kashiwara filtration of $G_{}^{\rm an}$ at $0$. By uniqueness it satisfies
\begin{eqnarray*}
V_\beta G_{}^{\rm an}&=&\CC\{\tau\}\otimes_{\CC[\tau]}V_\beta G
\end{eqnarray*}
Because $G$ is regular at $0$, there exists an isomorphism of $\CC\{\tau\}[\tau_{}^{-1}]$-connections
\begin{eqnarray}\label{eqnisom}
\lefcro \ooplus_{\beta\in\QQ}\gr_{\beta}^{V}G\rigcro_{}^{\rm an}\simeq G_{}^{\rm an}
\end{eqnarray}
which induces identity after graduation by the corresponding filtrations $V$ (see \eg. \cite[Prop\ptbl II.1.4]{Malgrange91} or also \cite[Lemma 6.2.6]{Bibi90}).

\begin{lemme}
For every $\beta\in\QQ$, $V_\beta G$ is a free $\CC[\tau]$-module of rank $\mu$ and defines a logarithmic connection with pole at $0$ on $U_0$.
\end{lemme}

\begin{proof}
As $V_\beta G$ has no $\CC[\tau]$-torsion, it is enough to prove that $V_\beta G$ has finite type over $\CC[\tau]$. This follows from the following two facts.

$\bullet$ The localization $\Ct\ootimes_{\CC[\tau]}V_\beta G$ has finite type over $\Ct$: this is due to (\ref{eqnresb}).

$\bullet$ The localization $\CC[\tau]_{(0)}^{}\otimes_{\CC[\tau]}V_\beta G$ has finite type over $\CC[\tau]_{(0)}^{}$: using faithful flatness of $\CC\{\tau\}$ over $\CC[\tau]_{(0)}^{}$ it is enough to consider the analytic localization, for which the assertion follows from the remark above, because $G$ is {\em regular} at $0$ and because it is clearly true for a $V$-graded module.
\end{proof}

\subsection*{The vector bundles $\ccG_{\beta,k}$}
Let us fix $\beta\in\QQ$ and $k\in \ZZ$. The locally free sheaf $\ccG_{\beta,k}$ is obtained by glueing $V_\beta G$ on $U_0$ and $G_k$ on $U_\infty$ using the isomorphisms on $U_0\cap U_\infty$
\begin{eqnarray*}
\Ct\ootimes_{\CC[\tau]}V_\beta G= G= \CC[\theta,\theta_{}^{-1}]\ootimes_{\CC[\theta]}G_k.
\end{eqnarray*}

The following is easy to prove.

\begin{lemme}\label{lemfacile}\mbox{ }
\begin{enumerate}
\item
We have $\ccG_{\beta,k}\otimes\cO(\ell)\simeq\ccG_{\beta+\ell,k}\simeq \ccG_{\beta,k+\ell}$ and $H^0(\PP^1,\ccG_{\beta,k})=V_\beta G\cap G_k$ where the intersection is taken in $G$. The isomorphism $H^0(\PP^1,\ccG_{\beta,k})\simeq H^0(\PP^1,\ccG_{\beta+k,0})$ is given by
\begin{eqnarray*}
\theta^k=\tau_{}^{-k}:V_\beta G\cap G_k\isom V_{\beta+k}G\cap G_0.
\end{eqnarray*}
\item
For every $\beta\in\QQ$ and $k\in\ZZ$, $V_\beta G\cap G_k$ is a finite dimensional vector space. For a fixed $\beta$ and for $k\ll0$ we have $V_\beta G\cap G_k=0$ and for $k\gg0$ we have $G_k=V_\beta G\cap G_k+G_{k-1}$. \hfill\qed
\end{enumerate}
\end{lemme}

In the following we will denote $\ccG_\beta\defin\ccG_{\beta,0}^{}$.

\subsection*{The spectral polynomial and the spectrum}
For $\beta\in\QQ$, let $$\nu_\beta=\dim V_\beta G\cap G_0/\lefpar V_\beta\cap G_{-1}+V_{<\beta}\cap G_0\rigpar.$$ The set of pairs $\{\beta,\nu_\beta\}$ for which $\nu_\beta\neq0$ is called the {\em spectrum} of $(G,G_0)$. The spectral polynomial of $(G,G_0)$ is
\begin{eqnarray*}
\SP_\psi(G,G_0;S)&\defin&\prod_{\beta\in\QQ}(S+\beta)_{}^{\nu_\beta}.
\end{eqnarray*}
If we put $\delta_\beta=\dim V_\beta G\cap G_0$, we have
\begin{eqnarray}
\nu_\beta=(\delta_\beta-\delta_{\beta-1})-(\delta_{<\beta}-\delta_{<\beta-1}).
\end{eqnarray}
In fact $V_\beta G$ induces a filtration on $G_0/G_{-1}$, namely \begin{eqnarray*}
V_\beta(G_0/G_{-1})&\defin& V_\beta G\cap G_0/V_\beta G\cap G_{-1},
\end{eqnarray*}
the dimension of which is $\delta_\beta-\delta_{\beta-1}$, and $\nu_\beta$ is the dimension of the $\beta$th graded piece. Hence the datum of the $\delta_\beta$'s is equivalent to the datum of the $\nu_\beta$'s.

\begin{lemme}
The degree of the spectral polynomial is equal to $\mu$.
\end{lemme}

\begin{proof}
Because we know that $\dim G_0/G_{-1}=\mu$, it is enough to show that
\begin{eqnarray*}
V_\beta (G_0/ G_{-1})=\lefacc
\begin{array}{ll}
0& \mbox{for $\beta\ll0$}\\
G_0/G_{-1}& \mbox{for $\beta\gg0$}
\end{array}
\rig
\end{eqnarray*}
which follows from lemma \ref{lemfacile}.
\end{proof}

Because $G$ is regular at $0$, the characteristic polynomial of the monodromy $T_0$ of the local system $\DRa G$ on $U_0^*$ is $\prod_{\alpha\in{}[0,1[}(T-\exp2i\pi\alpha)_{}^{\mu_\alpha}$ with $\mu_\alpha=\dim_\CC\gr_{\alpha}^{V}G$ (this follows from \cite[Prop\ptbl II.1.4]{Malgrange91}). Grading by $G_\bbullet$ we hence get
\begin{proposition}\label{proppolcar}
Let $\prod_\beta(S+\beta)_{}^{\nu_\beta}$ be the spectral polynomial of $(G,G_0)$. Then the polynomial $\prod_\beta(T-\exp2i\pi\beta)_{}^{\nu_\beta}$ is the characteristic polynomial of the monodromy $T_0$ of $\psi_\tau\DRa G$. \hfill\qed
\end{proposition}

\section{How to obtain a lattice by Fourier transform}\label{sec2}

Let $M$ be a regular holonomic $\CC[t]\langle\partial_t\rangle$-module (regularity at infinity is also assumed and, in fact, only regularity at infinity will be useful in prop\ptbl \ref{propfourier}) and put $M[\partial_{t}^{-1}]\defin\CC[\partial_t,\partial_{t}^{-1}]\otimes_{\CC[\partial _t]}^{}M$. It is also regular holonomic and moreover the (left) action of $\partial_t$ on $M[\partial_{t}^{-1}]$ is bijective, \ie. one can ``integrate in a unique way'' all elements of $M[\partial_{t}^{-1}]$.

We denote $\wh M$ the Fourier transform of $M$: it is equal to $M$ as a $\CC$-vector space, and is equipped with a left action of $\CC[\tau]\langle\partial _\tau\rangle$, where $\tau$ acts like $\partial _t$ and $\partial _\tau$ as $-t$ (see \eg. \cite{Malgrange91}).

Let $G$ be the Fourier transform of $M[\partial_{t}^{-1}]$. Then $G=\wh M[\tau_{}^{-1}]$ satisfies the assumptions made at the beginning of section \ref{sec1}. We will put as above $\theta=\tau_{}^{-1}=\partial_{t}^{-1}$ and will view $G$ as a $\CC[\theta,\theta_{}^{-1}]$-module with an action of $t=\theta^2\partial_\theta$.

We shall now show how to construct a lattice starting with a $\CC[t]$-module inside $M$. This will give a Fourier correspondence at the level of lattices.

\begin{proposition}\label{propfourier}
Let $M$ be a regular holonomic $\CC[t]\langle\partial_t\rangle$-module (regularity at infinity would be enough) and let $M_0$ be a finite type $\CC[t]$-submodule of $M$. Let $M'_0$ be the image of $M_0$ in $M[\partial_{t}^{-1}]$ by the natural morphism $M\rightarrow M[\partial_{t}^{-1}]$. Let $G_0\defin\sum_{i\geq 0}^{}\partial _{t}^{-i}M'_0=\CC[\theta]\cdot M'_0$ be the $\CC[\theta]$-module generated by $M'_0$. Then
\begin{enumerate}
\item
$G_0$ has finite type over $\CC[\theta]$;
\item
if $M_0$ generates $M$ over $\CC[t]\langle\partial_t\rangle$, then $G_0$ is a lattice in $G=\wwwh{M[\partial _{t}^{-1}]}$.
\end{enumerate}
\end{proposition}

\begin{proof}
Let us begin by $(1)$. First, consider the case $M=\CC[t]\langle\partial_t\rangle/(P)$ where $P=\sum_{i=0}^{d}\partial _{t}^{i}a_i(t)$, with $a_i\in\CC[t]$, $a_d\not\equiv0$ and (regularity at infinity) $\deg a_i<\deg a_d$ for $i<d$. The image $g$ of $1$ in $M[\partial_{t}^{-1}]$ thus satisfies $\partial_{t}^{-d}P\cdot g=0$, which can be written, putting $\theta=\partial_{t}^{-1}$, $\sum_{k=0}^{\mu}b_k(\theta)t^k\cdot g=0$, where $\mu=\deg a_d$, and with $b_i\in\CC[\theta]$ and $b_\mu$ is a nonzero constant (the leading coefficient of $a_d$). Hence the property is true for the $\CC[t]$-module $M_0$ generated by $1$.

Next, remark that if the property is true for some $M_0$ generating $M$, it is true for any $M_0$: one uses the fact that
$$
\CC[\theta]\cdot\lefpar M'_0+\partial_t M'_0\rigpar = G_0+\theta_{}^{-1}G_0=\theta_{}^{-1}G_0.
$$
Hence, $(1)$ is true for $M$ as above. Any $M$ comes in an exact sequence of $\CC[t]\langle\partial_t\rangle$-modules
$$
0\longrightarrow K\longrightarrow \CC[t]\langle\partial_t\rangle/(P)\longrightarrow M\longrightarrow 0
$$
with $P$ as above and $K$ a $\CC[t]$-torsion module: indeed, any holonomic $M$ is cyclic, so isomorphic to $\CC[t]\langle\partial_t\rangle/I$, for some left ideal $I$; take for $P$ an element of $I$ which has minimal degree with respect to $\partial_t$. The result being easy for torsion modules, it is then true for any $M$.

\medskip
$(2)$ is proved using the formula
\begin{eqnarray*}
\CC[\theta]\cdot\lefpar M'_0+\cdots+\partial_{t}^{k}M'_0\rigpar &=&\theta^{-k}G_0.
\end{eqnarray*}
\end{proof}

\begin{rem}\label{remrang}
For $M_0$ generating $M$, we have $\CC(t)\otimes_{\CC[t]}^{}M_0=\CC(t)\otimes_{\CC[t]}^{}M$ (and these are finite dimensional $\CC(t)$-vector spaces).
\end{rem}

Notice also that, if $M=M[\partial_{t}^{-1}]$, the dimension of this vector space is equal to the rank of $G$ as a $\CC[\tau,\tau_{}^{-1}]$-module: the rank of $G$ is equal to the sum of dimensions of vanishing cycles of $\pDRa M$ (see \eg. \cite[Cor\ptbl 8.3]{Brylinski86}); this is equal to the generic rank of $\pDRa M$ thanks to the fact that $\bH^j(\Afu,\pDRa(M))=0$ for all $j$ (this follows from the fact that $\partial _t:M\rightarrow M$ is an isomorphism and from the comparison theorem for $M$, which is regular holonomic).

\subsection*{Fourier transform and formal microlocalization}
We will adapt below the results of \cite[\T5]{Malgrange81} to any regular holonomic $\CC[t]\langle\partial _t\rangle$-module $M$. In \loccit. they are proved in the case where $M$ has only one singular point in $\Afu$, and this one is even not assumed to be regular (but the one at infinity is so); moreover the result is proved in a Gevrey context, not only in the formal one below, which will be enough for our purpose. Proposition \ref{propfourmicro} below is certainly ``well known to specialists'', but does not seem to exist in such a form in the literature.

Denote as above $\theta=\partial _{t}^{-1}$ and let $K=\CC\lcr\theta\rcr[\theta_{}^{-1}]$ be the ring of formal Laurent series in $\theta$. Put $t=\theta^2\partial _\theta$ and identifiy $K\otimes_{\CC[t]}^{}\CC[t]\langle\partial _t\rangle$ with $K\langle\partial _\theta\rangle$.

Let $\cE_{\afua}^{}$ be the sheaf of formal microdifferential operators on $T^*\Afua\setminus\mbox{zero section}$. We will consider it as a sheaf on $\Afua$ by restricting it to the section image($dt$). A local section of $\cE_{\afua}^{}$ is a formal Laurent series $\sum_{i\geq i_0}^{}a_i(t)\theta^i$ where $a_i$ are holomorphic functions defined on a fixed open set. We denote as usual $\cE_{\afua}^{}(0)$ the subring of sections of $\cE_{\afua}^{}$ with no pole at $\theta=0$.

Let $M$ be a $\CC[t]\langle\partial _t\rangle$-module and put $G=M[\partial_{t}^{-1}]$ viewed as a $\CC[\tau]\langle\partial_\tau\rangle$-module, or as $\CC[\theta]\langle\partial_\theta\rangle$-module. Put also $\cM=\cO_{\afua}^{}\otimes_{\CC[t]}^{}M$ and let $\muM=\cE_{\afua}^{}\otimes_{\cD_{\afua}^{}}^{}\cM$ be its formal microlocalization. As a sheaf on $\Afua$, it is supported on the singular set of $M$, since $\muM=0$ when $\cM$ is $\cO_{\afua}^{}$-locally free. At a singular point $c$, the action of $t$ on the germ $\muM_c$ can be written as $e_{}^{c/\theta}(t-c)e_{}^{-c/\theta}$ if $\cE_c$ is viewed as a subring of $\CC\{t-c\}\lcr\theta\rcr$.

Let $M_0$ be a finite type $\CC[t]$-submodule of $M$ generating $M$ over $\CC[t]\langle\partial_t\rangle$ and put $\cM_0=\cO_{\afua}^{}\otimes_{\CC[t]}^{}M_0$. Then $\muM_0\defin{\rm image}\lefcro \cE(0)\otimes_{\cO}^{}\cM_0\rightarrow \muM\rigcro$ is a lattice in $\muM$. From the preparation theorem for $\cE(0)$-modules (see \eg. \loccit.) follows that the germ $\muM_{0,c}$ at each singular point $c$ of $\muM_0$ is a free $\CC\lcr\theta\rcr$-module and that $\muM_c=K\otimes_{\CC\lcr\theta\rcr}^{}\muM_{0,c}$ is a finite $K$-vector space with an action of $\partial _\theta$.

On the other hand, let $G_0$ be as above and denote $G_{0}^{\wh{\hphantom{0}}}=\CC\lcr\theta\rcr\otimes_{\CC[\theta]}^{}G_0$. As $G_{0}^{\wh{\hphantom{0}}}$ is stable by $t=\theta^2\partial _\theta$, it follows from the method of the Turrittin theorem on the splitting of formal connections (see \eg. \cite[Chap\ptbl III]{Malgrange91}) that the pair $(G_{}^{\wh{\hphantom{0}}},G_{0}^{\wh{\hphantom{0}}})$ splits as a direct sum indexed by the singular points of $M$:
\begin{eqnarray*}
(G_{}^{\wh{\hphantom{0}}},G_{0}^{\wh{\hphantom{0}}})&\simeq&\ooplus_c (G_{c}^{\wh{\hphantom{0}}}\otimes e_{}^{-c/\theta},G_{0,c}^{\wh{\hphantom{0}}}\otimes e_{}^{-c/\theta})
\end{eqnarray*}
where each $G_{c}^{\wh{\hphantom{0}}}$ is a regular $K$-connection and $G_{0,c}^{\wh{\hphantom{0}}}$ is a lattice in it; moreover, $\otimes e_{}^{-c/\theta}$ denotes the twist of the $\partial_\theta$ (or the $\theta^2\partial_\theta$) action.

\begin{proposition}\label{propfourmicro}
The composed $\CC\lcr\theta\rcr$-linear map
$$
G_{}^{\wh{\hphantom{0}}}\defin K\ootimes _{\CC[\theta_{}^{-1}]}^{}M\longrightarrow \Gamma\lefpar \Afua,K\ootimes_{\CC[\partial_t]}^{}\cM\rigpar\longrightarrow \Gamma\lefpar \Afua,\muM\rigpar 
$$
is an isomorphism compatible with the action of $t$, which identifies $G_{0}^{\wh{\hphantom{0}}}$ with $\Gamma\lefpar \Afua,\muM_0\rigpar$.
\end{proposition}

\begin{proof}
Remark first that
$$
K\ootimes _{\CC[\theta_{}^{-1}]}^{}M=K\ootimes _{\CC[\theta,\theta_{}^{-1}]}^{}\lefpar \CC[\theta,\theta_{}^{-1}]\ootimes _{\CC[\theta_{}^{-1}]}^{}M\rigpar =K\ootimes _{\CC[\theta,\theta_{}^{-1}]}^{}G=\CC\lcr\theta\rcr\ootimes_{\CC[\theta]}^{}G.
$$
It is known that $G_{}^{\wh{\hphantom{0}}}$ and $\Gamma\lefpar \Afua,\muM\rigpar$ are $K$-vector spaces of the same dimension, namely the sum of dimensions of vanishing cycles of $\DRa(\cM)$ at its singularities (see remark \ref{remrang} above). Moreover, the map is clearly compatible with the action of $t$. Last, the map clearly sends $G_{0}^{\wh{\hphantom{0}}}$ into $\Gamma(\Afua,\muM_0)$. As there is no nonzero morphism compatible with the $t$ action between $G_{c}^{\wh{\hphantom{0}}}\otimes e_{}^{-c/\theta}$ and $\muM_{c'}$ for $c\neq c'$, it is enough to show that for each $c$ the map $G_{0}^{\wh{\hphantom{0}}}\rightarrow \muM_{0,c}$ is onto, or, equivalently, due to Nakayama's lemma, that $G_{0}^{\wh{\hphantom{0}}}\rightarrow \muM_{0,c}/\theta\muM_{0,c}$ is onto.

Remark that, since $\cM_{0,c}^{}=\cO_{c}^{}\otimes_{\CC[t]}^{}M_0$, there exists $m_1,\ldots ,m_p\in M_0$ such that any $m\in\cM_{0,c}^{}$ can be written $\sum \varphi_im_i$ with $\varphi_i\in\cO_c$. There exists $d\in\NN$ and $a_d(t)\in\CC[t]$ such that $\partial_{t}^{d}a_d(t)\cdot m_i\in\sum_{k=1}^{d}\partial _{t}^{d-k}M_0$ for all $i=1,\ldots ,p$. If $n$ is the order of vanishing of $a_d$ at $c$, we conclude that $(t-c)^n1\otimes m_i\in\theta\muM_{0,c}$ for all $i$. As we have $m=m'+(t-c)^n\sum\psi_im_i$ with $m'\in M_0$ and $\psi_i\in\cO_c$, we conclude that $1\otimes m\equiv1\otimes m' \mod \theta\muM_0$ in $\muM_0$.
\end{proof}

\section{Spectrum and duality}
The free module $G^*\defin \Hom_{\ct}(G,\Ct)$ is naturally equipped with a structure of a meromorphic connection: put, for $\varphi\in G^*$,
\begin{eqnarray}\label{eqnhom}
(\dtau \varphi)(g)= \dtau (\varphi(g))-\varphi(\dtau g).
\end{eqnarray}

\begin{lemme}\label{lemdualspectre}
For $\beta\in\QQ$, we have $V_\beta(G^*)=\Hom_{\CC[\tau]}(V_{<-\beta+1}G,\CC[\tau])$.
\end{lemme}

\begin{proof}
Put $\ccV_\beta G^*=\Hom_{\CC[\tau]}(V_{<-\beta+1}G,\CC[\tau])$. This is a free $\CC[\tau]$-module of rank $\mu$ naturally contained in $G^*$ and stable under the action of $\tau\dtau$. Moreover we have $\tau^k\ccV_\beta G^* = \ccV_{\beta-k}G^*$. Let $N$ such that for all $\delta$, $(\tau\dtau+\delta)^NV_{\delta}G\subset V_{<\delta}G$. Let $\varphi\in\ccV_\beta G^*$. We will now show that $(\tau\partial_\tau+\beta)^N\varphi\in\ccV_{<\beta}^{}G^*=\Hom_{\CC[\tau]}(V_{-\beta+1}G,\CC[\tau])$. Remark that for $g\in V_{-\beta+1}G$ we have
\begin{eqnarray*}
\lefcro (\tau\partial_\tau+\beta)\varphi\rigcro(g)&=&\partial_\tau \varphi(\tau g)-\varphi\lefpar (\tau\partial_\tau-\beta+1)g\rigpar 
\end{eqnarray*}
and $\partial_\tau \varphi(\tau g)\in\CC[\tau]$ since $\tau g\in V_{-\beta}G\subset V_{<-\beta+1}G$. By induction we see that
\begin{eqnarray*}
\lefcro (\tau\partial_\tau+\beta)^N\varphi\rigcro(g)&=&h+(-1)^N\varphi\lefpar (\tau\partial_\tau-\beta+1)^Ng\rigpar 
\end{eqnarray*}
with $h\in\CC[\tau]$. The assertion follows from the fact that $(\tau\partial_\tau-\beta+1)^Ng\in V_{<-\beta+1}G$.

Let us verify that the filtration $\ccV_\bbullet G^*$ is good relatively to $V_\bbullet\CC[\tau]\langle\partial_\tau\rangle$. We have to verify that for each $\beta$ we have $\ccV_{\beta+k}^{}G^*=\partial_\tau\ccV_{\beta+k-1}^{}G^*+\ccV_{\beta+k-1}^{}G^*$ for all $k\geq k_0$. This is equivalent to
\begin{eqnarray*}
\ccV_{\beta+k-1}^{}G^*&=&\tau\partial_\tau \ccV_{\beta+k-1}^{}G^*+\ccV_{\beta+k-2}^{}G^*
\end{eqnarray*}
and is true as soon as $\beta+k-1>0$ as a consequence of the previous result.

Now $\ccV_\bbullet G^*$ satisfies all the characteristic properties of the Malgrange-Kashiwara filtration $V_\bbullet G^*$, hence is equal to it.
\end{proof}

If $G_0$ is a lattice in $G$, then $G^*_0\defin\Hom_{\CC[\theta]}(G_0,\CC[\theta])$ can be identified with the set of $\varphi\in G^*$ verifying $\varphi(G_0)\subset \CC[\theta]$ or also $\varphi(G_k)\subset \theta_{}^{-k}\CC[\theta]\subset\CC[\theta,\theta_{}^{-1}]$ for all $k\in\ZZ$. It follows easily from (\ref{eqnhom}) that $G^*_0$ is a lattice in $G^*$. From the previous lemma we conclude that $\ccG_{\beta}^{*}$ is equal to the locally free sheaf dual to $\ccG_{<-\beta+1}$. 

\begin{proposition}\label{propdualspectre}
If $\prod_{\beta\in\QQ}(S+\beta)_{}^{\nu_\beta}$ is the spectral polynomial of $(G,G_0)$, then the polynomial $\prod_{\beta\in\QQ}(S-\beta)_{}^{\nu_\beta}$ is the one of $(G^*,G^*_0)$; in other words we have $\nu_\beta(G^*,G_0^*)=\nu_{-\beta}^{}(G,G_0)$.
\end{proposition}

\begin{proof}
Let us keep notation as above. Serre duality and the fact that $\ccG_\gamma^*=(\ccG_{<-\gamma+1}^{})^*$ show that
\begin{eqnarray*}
h^1(\PP^1,\ccG_{<-\gamma+1})=h^0(\PP^1,\ccG_{\gamma-2}^{*})& \mbox{and} &h^1(\PP^1,\ccG_{-\gamma+1})=h^0(\PP^1,\ccG_{<\gamma-2}^{*})
\end{eqnarray*}
Conclude by considering the exact sequences
$$
0\rightarrow H^0(\PP^1,\ccG_{<-\gamma})\rightarrow H^0(\PP^1,\ccG_{-\gamma})\rightarrow \gr_{-\gamma}^{V}G\rightarrow H^1(\PP^1,\ccG_{<-\gamma})\rightarrow H^1(\PP^1,\ccG_{-\gamma})\rightarrow 0
$$
for $\gamma=\beta+1,\beta$ and the fact that $\dim \gr_{-\beta}^{V}G=\dim \gr_{-\beta-1}^{V}G$.
\end{proof}

\begin{corollaire}\label{corsymspectre}
If there exists an isomorphism $G^*\isom G$ inducing $G_0^*\simeq G_w$ for some $w\in\ZZ$, then the spectrum of $(G,G_0)$ is symmetric relatively to $w/2$, \ie. $\nu_\beta=\nu_{w-\beta}$ for all $\beta\in\QQ$. \hfill\qed
\end{corollaire}

\begin{rem}\label{remsesqui}
Consider the involution $\CC[\tau]\langle\partial_\tau\rangle\isom \CC[\tau]\langle\partial_\tau\rangle$ given by $\tau\mapsto -\tau$, $\partial_\tau\mapsto -\partial_\tau$ and denote $\ov G$ the module $G$ where the action of operators is composed with this involution. Because the action of $\tau\partial_\tau$ is unchanged, we have $V_\alpha\ov G=\ov{V_\alpha G}$. Moreover, if $G_0$ is a lattice in $G$, then $\ov{G_0}$, which is equal to $G_0$ on which $\theta$ acts as $-\theta$ and $t$ as $-t$, is also a lattice in $\ov G$. In particular the spectrum of $(G,G_0)$ is equal to the one of $(\ov G,\ov{G_0})$. If we have an isomorphism $(\ov{G^*},\ov {G^*_0})\isom (G,G_w)$, we can apply corollary \ref{corsymspectre}.
\end{rem}

\begin{rem}\label{remdualspectre}
Proposition \ref{propdualspectre} can also be obtained in the following more concrete way. From lemma \ref{lemdualspectre} we deduce that the natural pairing $V_\beta G^*\times V_{-\beta+1}^{}G\rightarrow \CC[\tau,\tau_{}^{-1}]$ has image in $\tau_{}^{-1}\CC[\tau]$ and the pairing with values in $\CC$ obtained by composing the previous one with the residue at $\tau=0$ induces a perfect pairing $\gr_{\beta}^{V}G^*\times \gr_{-\beta+1}^{V}G\rightarrow \CC$.

If for any $k\in\ZZ$ we denote $G_k^*=\tau^kG_0^*$, we get in the same way an isomorphism
\begin{eqnarray*}
\gr_{k}^{G^*}\gr_{\beta}^{V}G^*&\simeq&\Hom_\CC(\gr_{-k-1}^{G}\gr_{-\beta+1}^{V}G,\CC)
\end{eqnarray*}
which implies \ref{propdualspectre}. We will have to use this isomorphism for $\beta=\alpha\in{} ]0,1[$, so that $-\alpha+1$ is also in $]0,1[$; for $\beta=0$ we will use the composition of this isomorphism by multiplication by $\tau$ to get
\begin{eqnarray*}
\gr_{k}^{G^*}\gr_{0}^{V}G^*&\simeq&\Hom_\CC(\gr_{-k}^{G}\gr_{0}^{V}G,\CC).
\end{eqnarray*}
\end{rem}

\subsection*{Behaviour with respect to tensor product}
The following property is useful in order to prove a Thom-Sebastiani type theorem for the spectrum. Let $(G',G'_0)$ and $(G'',G''_0)$ be as in section \ref{sec1} and put $G=G'\otimes_{\CC[\theta,\theta_{}^{-1}]}^{}G''$. Then $G$ is regular at $\tau=0$ and $G_0\defin G'_0\otimes_{\CC[\theta]}^{}G''_0$ is a lattice in it.

\begin{proposition}\label{proptensor}
Assume that we are given isomorphisms $$(\ov{{G'}^*},\ov {{G'}^*_0})\isom (G',G'_{w'})\quad\mbox{and}\quad (\ov{{G''}^*},\ov {{G''}^*_0})\isom (G'',G''_{w''})$$ for some $w',w''\in\ZZ$. Then the spectrum of $(G,G_0)$ is given by
\begin{eqnarray*}
\nu_\beta(G,G_0)&=&\sum_{\beta'+\beta''=\beta}^{}\nu_{\beta'}^{}(G',G'_0)\nu_{\beta''}^{}(G'',G''_0).
\end{eqnarray*}
\end{proposition}

\begin{proof}
The proof follows the one of \cite{Varchenko82} or \cite{S-S85}. Remark first that we have an isomorphism $(\ov{G^*},\ov {G^*_0})\isom (G,G_w)$ with $w=w'+w''$ by taking the tensor product of the one for $G'$ and the one for $G''$, as $G'_{w'}\otimes G''_{w''}=G_w$. As $G',G'',G$ are regular at $\tau=0$ we have (by first twisting with $\CC\{\tau\}$ and using the structure of regular connections)
\begin{eqnarray*}
\gr_\beta^VG&=&\ooplus_{\beta'\in[0,1[}^{}\lefpar \gr_{\beta'}^{V}G'\ootimes_\CC\gr_{\beta-\beta'}^{V}G''\rigpar.
\end{eqnarray*}
Consider on the LHS the filtration induced by $\tau^kG_0$ and on the RHS the tensor product of the corresponding filtrations on $\gr_{\beta'}^{V}G'$ and $\gr_{\beta-\beta'}^{V}G''$. The former contains the latter, so in particular we have
\begin{eqnarray*}
\sum_{k\geq 0}^{}\nu_{\beta-k}^{}=\dim G_0\gr_\beta^VG&\geq &\sum_{\beta'\in[0,1[}^{}\sum_{k\geq 0}^{}\sum_{i+j=k}^{}\nu'_{\beta'-i}\nu''_{\beta-\beta'-j}.
\end{eqnarray*}
Put $(\nu'\star\nu'')_\beta=\sum_{\beta'+\beta''=\beta}^{}\nu'_{\beta'}\nu''_{\beta''}$. Then we conclude that $\sum_{k\geq 0}^{}\nu_{\beta-k}^{}\geq \sum_{k\geq 0}^{}(\nu'\star\nu'')_{\beta-k}^{}$ and thus $\sum_{\gamma\leq \beta}^{}\nu_{\gamma}^{}\geq \sum_{\gamma\leq \beta}^{}(\nu'\star\nu'')_{\gamma}^{}$ for any $\beta$. As both terms are equal for $\beta\ll0$ or $\beta\gg0$ and as $\nu$ and $\nu'\star\nu''$ are both symmetric with respect to $w/2$, we conclude that $\nu=\nu'\star\nu''$.
\end{proof}

\subsection*{A microdifferential criterion for the symmetry of the spectrum}
We keep notation of \T\ref{sec2}: the module $M$ is $\CC[t]\langle \partial_t\rangle$-holonomic and regular even at infinity and we assume that $M$ comes equipped with a $\CC[t]$-module $M_0$ of finite type generating $M$.

Let $DM$ be the $\CC[t]\langle \partial_t\rangle$-module dual to $M$, \ie. the left module associated with the right module $$\Ext_{\CC[t]\langle \partial_t\rangle}^{1}(M,\CC[t]\langle \partial_t\rangle).$$ We also denote $D(\muM)$ the dual of the microdifferential system attached to $M$, which is identified with $(D\cM)^{\!{}^\mu}$. As $\muM$ is endowed with a good filtration $\muM_\bbullet=\cE(\bbullet)\cdot\muM_0$, the dual $D(\muM)$ is equipped with a natural filtration $D(\muM)_\bbullet$ defined using any strictly filtered resolution of $(\muM,\muM_\bbullet)$ by free $(\cE,\cE(\bbullet))$-modules.

Assume that we are given a morphism $DM\rightarrow M$ such that the kernel and the cokernel are free $\CC[t]$-modules of finite rank. Taking Fourier transforms and localizing with respect to $\tau$ one gets an isomorphism $\wh{DM}[\tau_{}^{-1}]\isom\wh M[\tau_{}^{-1}]=G$. But we have the following relation between duality and Fourier transform (see \cite[lemme V.3.6]{Malgrange91}):
\begin{eqnarray*}
D\wh M=\ov {\wh{DM}}
\end{eqnarray*}
where $D$ denotes the duality as $\CC[t]\langle\partial _t\rangle$ or $\CC[\tau]\langle\partial_\tau\rangle$-module. We will identify below $(D\wh M)[\tau_{}^{-1}]$ with $G^*$ defined above. Hence we get an isomorphism $\ov{G^*}\isom G$.

\begin{lemme}[{\cite[\T2.7]{MSaito89}}]\label{lemdualsaito}
We have a canonical isomorphism $$G^*=(DG)[\tau_{}^{-1}]=(D\wh M)[\tau_{}^{-1}].$$
\end{lemme}

\begin{proof}
The second equality is clear as the kernel and cokernel of $\wh M\rightarrow G$ are supported at $\tau=0$. We have
\begin{eqnarray*}
\CC[\tau,\tau_{}^{-1}]\ootimes_{\CC[\tau]}^{}DG&=&\Ext_{\CC[\tau,\tau_{}^{-1}]\langle\partial_\tau\rangle}^{1}(G,\CC[\tau,\tau_{}^{-1}]\langle\partial_\tau\rangle)^g
\end{eqnarray*}
where  $N^g$ means the left module associated with the right module $N$; consider then the following resolution of $G$ as a left $\CC[\tau,\tau_{}^{-1}]\langle\partial_\tau\rangle$-module:
\begin{eqnarray*}
&0\rightarrow \CC[\tau,\tau_{}^{-1}]\langle\partial_\tau\rangle\ootimes_{\CC[\tau,\tau_{}^{-1}]}^{} G\MRE{\partial_\tau\otimes 1-1\otimes \partial_\tau}\CC[\tau,\tau_{}^{-1}]\langle\partial_\tau\rangle\ootimes_{\CC[\tau,\tau_{}^{-1}]}^{} G\longrightarrow G\rightarrow 0&
\end{eqnarray*}
with $\CC[\tau,\tau_{}^{-1}]\langle\partial_\tau\rangle\ootimes_{\CC[\tau,\tau_{}^{-1}]}^{} G= \CC\langle\partial_\tau\rangle\ootimes_{\CC}^{} G$, and identify
$$
\Hom_{\CC[\tau,\tau_{}^{-1}]\langle\partial_\tau\rangle}^{}\lefpar \CC\langle\partial_\tau\rangle\ootimes_{\CC}^{} G,\CC[\tau,\tau_{}^{-1}]\langle\partial_\tau\rangle\rigpar ^g=(G^*\otimes_\CC\CC\langle\partial_\tau\rangle)^g=\CC\langle\partial_\tau\rangle\otimes_\CC G^*
$$
to get the assertion.
\end{proof}

On the other hand, the morphism $DM\rightarrow M$ induces an isomorphism $D(\muM)\rightarrow \muM$. When the microdifferential property below holds, the spectrum of $(G,G_0)$ is symmetric with respect to $w/2$.

\begin{proposition}\label{propresolstricte}
In the previous situation, assume moreover that there exists $w\in\ZZ$ such that this isomorphism sends $D(\muM)_0$ onto $\theta_{}^{-w}\muM_0=\muM_w$. Then the isomorphism $\ov{G^*}\rightarrow G$ sends $\ov{G^*_0}$ onto $G_w$.
\end{proposition}

\begin{proof}
As the image of $\ov{G^*_0}$ is contained in some $G_k$ and coincides with $G$ after localizing with respect to $\tau$, it is enough to show that $\CC\lcr\theta\rcr\otimes_{\CC[\theta]}^{}\ov{G^*_0}$ is sent onto $\CC\lcr\theta\rcr\otimes_{\CC[\theta]}^{}G_w$. According to proposition \ref{propfourmicro}, the problem is reduced to identifying, for each singular point $c$, $(\Hom_K(\muM_c,K),\Hom_{\CC\lcr\theta\rcr}^{}(\muM_{0,c}\CC\lcr\theta\rcr))$ with $(D(\muM_c),D(\muM_c)_0)$. This follows also from \cite[\T2.7]{MSaito89}: recall that $\muM_c$ is a finite dimensional $K$-vector space, so one has a resolution of $\muM_c$ as an $\cE_c$-module
$$
0\longrightarrow \cE_c\ootimes_K\muM_c\MRE{(t-c)\otimes 1-1\otimes (t-c)}\cE_c\ootimes_K\muM_c\longrightarrow \muM_c\longrightarrow 0
$$
which induces a resolution
$$
0\longrightarrow \cE_c(0)\ootimes_{\CC\lcr\theta\rcr}\muM_{0,c}\MRE{(t-c)\otimes 1-1\otimes (t-c)}\cE_c(0)\ootimes_{\CC\lcr\theta\rcr}\muM_{0,c} \longrightarrow \muM_{0,c}\longrightarrow 0.
$$
In fact it is easier to show first that the corresponding sequences are exact when $\cE_c,\cE_c(0)$ are replaced with $K\langle\partial _\theta\rangle,\CC\lcr\theta\rcr\langle\theta^2\partial _\theta\rangle$, and then to use flatness of the former rings over the latter ones respectively.
\end{proof}

\section{The Riemann-Hilbert-Birkhoff problem}

We keep notation of the previous section. In order to simplify notation, it will be convenient to write a basis as a column vector; the matrices considered below are the transposed matrices of the usual matrices.

\begin{proposition}\label{propequiv}
The following properties are equivalent:
\begin{enumerate}
\item
there exists a basis $\bmg$ of $G_0$ over $\CC[\theta]$ for which the matrix of $t=\theta^2\dth$ is equal to $A_0+\theta A_1$ where $A_0$ and $A_1$ are constant matrices;
\item
there exists a free $\CC[\tau]$-module $H^0\subset G$ stable under the action of $\tau\dtau$ generating $G$ over $\Ct$ and such that $(H^0\cap G_0)\oplus G_{-1}=G_0$;
\item
there exists a free $\cO_{\PP^1}$-module $\cG_0$ inside (the $\cO_{\PP^1}$-module associated with) $G$ such that the restriction of $\cG_0$ on $U_\infty$ is equal to $G_0$ and the restriction to $U_0$ is a logarithmic connection with pole at $0$.
\end{enumerate}
\end{proposition}

\begin{rems*}\mbox{ }
\begin{enumerate}
\item
If such a basis exists, then $A_0$ is the transposed matrix of the action of $t$ on $G_0/G_{-1}$ and $-A_1$ is the transposed matrix of the residue at $\tau=0$ of the logarithmic connection $H^0$; moreover the trace of $A_1$ does not depend on the choice of such a basis: in fact, given two bases $\bmg$ and $\gammag$, one can assume that they induce the same basis on $G_0/G_{-1}$ up to a constant base change on $\gammag$, so the matrices $A_0$ are the same; there exists $P\in \GL_\mu(\CC[\theta])$ such that
\begin{eqnarray*}
A_0+\theta A'_1&=&\theta^2\dth (P)\cdot P_{}^{-1}+ P(A_0+\theta A_1)P_{}^{-1};
\end{eqnarray*}
as $\det P$ is a nonzero constant, this implies $\tr A'_1=\tr A_1$.
 \item
In general, one knows the existence of such a basis in the following cases:
\begin{itemize}
\item
the monodromy matrix of $G$ (around $0$ or $\infty$) has $\mu$ distinct eigenvalues (see for instance \cite{Malgrange83c});
\item
$G_0$ is semi-simple as a $\CC[\theta]$-module equipped with the action of $t$ (\cf. \cite{Bolibruch94}).
\end{itemize}
\end{enumerate}
\end{rems*}

\begin{proof}
$(1)\implique(2)$: one puts $H^0=\CC[\tau]\cdot \bmg$. Then $H^0$ is stable under the action of $\tau\dtau$ and the matrix of $\tau\dtau$ in the basis $\bmg$ is $-(\tau A_0+A_1)$. In particular $A_1$ is the residue of the logarithmic connection on $H^0$. One has $\Ct\otimes H^0=G$ and $H^0\cap G_0=\CC\cdot \bmg$.

$(2)\implique(3)$: because of the isomorphisms $H^0[\tau_{}^{-1}]=G=G_0[\theta_{}^{-1}]$, one may construct a locally free $\cO_{\PP^1}$-module $\cG_0$ such that $\cG_{0|U_0}=H^0$ and $\cG_{0|U_\infty}=G_0$. We have $\Gamma(\PP^1,\cG_0)=H^0\cap G_0$ and the assumption means that the restriction morphism $\Gamma(\PP^1,\cG_0)\rightarrow i_\infty^*\cG_0$ is an isomorphism. From the Birkhoff-Grothendieck theorem one deduces easily that $\cG_0$ is trivial.

$(3)\implique(1)$: let $\bmg$ be a basis of $\cG_0$. Because $G_0$ is stable under the action of $t$, the matrix of $t$ in this basis has elements in $\CC[\theta]$. Because $H^0\defin \cG_{0|U_0}$ is stable under the action of $\tau\dtau=-\theta\dth$, one concludes that this matrix has degree at most one in $\theta$.
\end{proof}

Assume that $(G,G_0)$ is obtained from $(M,M_0)$ as in \T\ref{sec2}, where $M=M[\partial_{t}^{-1}]$ is regular holonomic, generated by $M_0$, and $M_0$ is $\CC[t]$-free and stable under $\partial_{t}^{-1}$. The {\em partial Riemann-Hilbert problem} for $M_0$ consists in finding a $\CC[t]$-basis of $M_0$ for which the connection matrix has at most a logarithmic pole at infinity.

\begin{corollaire}\label{corRH}
Assume that $M=M[\partial_{t}^{-1}]$ is regular even at infinity, that $M_0$ has finite type over $\CC[t]$, is stable under $\partial_{t}^{-1}$ and generates $M$. Then, if the R-H-B problem has a solution for $G_0$ with $A_1+k\id$ invertible for all $k\in\NN$, the $\CC[t]$-module $M_0$ is free and the partial Riemann-Hilbert problem has a solution for $M_0$.
\end{corollaire}

\begin{proof}
Indeed, let $\bmg$ be a $\CC[\theta]$-basis of $G_0$ satisfying property (1) of proposition \ref{propequiv}. Then we have, for all $k\geq 1$,
\begin{eqnarray*}
\theta^k\bmg&=&\prod_{\ell=0}^{k-1}\lefcro \lefpar A_{1}+\ell\id\rigpar ^{-1} (t\id-A_0)\rigcro\cdot\bmg, 
\end{eqnarray*}
hence $\bmg$ generates $M_0$ as a $\CC[t]$-module. Since $\rk G_0=\rk M_0$ (see remark \ref{remrang}), $M_0$ is $\CC[t]$-free and $\bmg$ is a $\CC[t]$-basis of $M_0$. Now the relation $t\bmg=(A_0+\theta A_1)\bmg$ can also be written $\partial _t\bmg=(A_0-t\id)_{}^{-1}(A_1-\id)\bmg$ and the connection matrix in this basis has rational coefficients with a pole of order at most one at infinity.
\end{proof}

\section{Good bases}
The content of this section is an adaptation of \cite[\T3]{MSaito89}.

\begin{definition}\label{defbb}
We shall say that a basis $\bmg$ of $G_0$ is a {\em good basis} if
\begin{enumerate}
\item
the matrix of $t=\theta^2\dth$ in this basis is $A_0+\theta A_1$, where $A_0$ and $A_1$ are constant $\mu\times \mu$ matrices;
\item
the spectrum of $A_1$ is equal to the spectrum of $(G,G_0)$.
\end{enumerate}
If moreover $A_1$ is semisimple, we shall say that $\bmg$ is a {\em very good basis}.
\end{definition}

\begin{proposition}\label{propbbases}
Assume that for every $\alpha\in{} [0,1[$ there exists a decreasing filtration $(H_{\alpha}^{k})_{k\in\ZZ}$ of $\gr_{\alpha}^{V}G$ by subspaces satisfying the following properties:
\begin{enumerate}
\item
for every $\alpha\in{} [0,1[$ and every $k\in\ZZ$, the subspace $H_{\alpha}^{k}$ is stable under the action of the nilpotent endomorphism $N$ of $\gr_{\alpha}^{V}G$ induced by $\tau\dtau+\alpha$;
\item (transversality condition)
the decreasing filtration $H_{\alpha}^{\bbullet}$ is opposite to the increasing filtration on $\gr_{\alpha}^{V}G$ induced by $(G_k)_{k\in\ZZ}$, \ie.
\begin{eqnarray*}
\gr_{H_\alpha}^{\ell}\gr_{k}^{G}\gr_{\alpha}^{V}G=0\quad\mbox{for }k\neq\ell.
\end{eqnarray*}
\end{enumerate}
Then one can construct a good basis of $G_0$. Moreover, this basis is very good if and only if for every $\alpha\in{} [0,1[$ the filtration $H_{\alpha}^{\bbullet}$ satisfies $NH_{\alpha}^{k}\subset H_{\alpha}^{k+1}$ for all $k\in\ZZ$.
\end{proposition}

\begin{proof}
Let $E=V_1G/V_{0}G$. This is a rank $\mu$ vector space equipped with an action of $\tau\dtau$. Moreover $E=\oplus_{\alpha\in{}[0,1[}E_\alpha$ with $E_\alpha=\cup_n \ker(\tau\dtau+\alpha)^n$ and by construction, the filtration induced by $V_\alpha G$ on $E$, namely $V_\alpha G/V_{0}G$ for $\alpha\in{}[0,1[$, is equal to $\oplus_{0<\alpha'\leq \alpha} E_{\alpha'}$. We shall denote it $V_\alpha E$. We have also $\gr_{\alpha}^{V}E=\gr_{\alpha}^{V}G$ for $\alpha\in{}[0,1[$.

The filtration $G_k$ induces a filtration $G_kE=G_k\cap V_1G/G_k\cap V_{0}G$ in the same way. The filtration $G_kE$ induces on $\gr_{\alpha}^{V}E=\gr_{\alpha}^{V}G$ the filtration $G_k\gr_{\alpha}^{V}G\defin G_k\cap V_\alpha G/G_k\cap V_{<\alpha}G$.

\begin{lemme}\label{lemrelev}
Under the assumption of proposition {\rm(\ref{propbbases})}, there exists a decreasing filtration $(L^k)_{k\in\ZZ}$ of $E$ which satisfies
\begin{enumerate}
\item
$L^k$ is stable under the action of $\tau\dtau$,
\item
$L_{}^{\bbullet}$ is opposite to $G_{\bbullet}$.
\end{enumerate}
\end{lemme}

\begin{proof}
We shall show (by induction on the number of $\alpha$ such that $E_\alpha\neq0$) that $L^k=\oplus_{\alpha}H_{\alpha}^{k}$ satisfies the desired properties. Put $\alpha_0=\max\{\alpha\in{}[0,1[{}\mid E_{\alpha}\neq0\}$, $E'=\oplus_{0\leq \alpha<\alpha_0}E_\alpha$ and $E''=E/E'=\gr_{\alpha_0}^{V}G$. We have to show that for every $k\in\ZZ$ we have
\begin{eqnarray*}
L^k\cap G_{k-1}E=0\quad\mbox{and}\quad L^k+G_{k-1}E=E.
\end{eqnarray*}
First, the induction hypothesis implies that $L^k\cap G_{k-1}E\cap E'=0$. Hence the map $L^k\cap G_{k-1}E\rightarrow E''$ is injective. But its image is contained in $H_{\alpha_0}^{k}\cap G_{k-1}\gr_{\alpha_0}^{V}G$, which is equal to $0$ by transversality assumption. The other equality is proved in the same way.
\end{proof}

In order to obtain information on $G$ from information on $E$, we shall introduce the Rees module of $G$ associated with the filtration $VG$. So let $u$ be a new variable and consider
\begin{eqnarray*}
\cR_{V}G_0\defin\ooplus_{k\in\ZZ}(V_kG\cap G_0)u^k\subset G_0\ootimes_\CC\CC[u,u_{}^{-1}].
\end{eqnarray*}
Denote $w$ the action of $\theta u$ on $G\ootimes_\CC\CC[u,u_{}^{-1}]$. Then $\cR_{V}G_0$ is naturally a $\CC[u,w]$-module because multiplication by $\theta$ induces an isomorphism
\begin{eqnarray*}
V_kG\cap G_0\isom V_{k+1}G\cap G_{-1}\subset V_{k+1}\cap G_0.
\end{eqnarray*}
Put $\partial_w=u_{}^{-1}\dth$ and $w^2\partial_w=u\theta^2\dth$. Then, because $\theta^2\dth$ acts as $-\dtau$ on $V_kG$, it induces a map
\begin{eqnarray*}
\theta^2\dth:V_kG\cap G_0\longrightarrow V_{k+1}G\cap G_0
\end{eqnarray*}
and one deduces an action of $w^2\partial_w$ on $\cR_{V}G_0$.

\begin{lemme}
$\cR_{V}G_0$ is a free $\CC[u,w]$-module of rank $\mu$. Its restriction at $u=u_0\neq0$ (with action of $w^2\partial_w$) is isomorphic to $G_0$ (with the action of $\theta^2\dth$) and its restriction at $u=0$ is equal to $\ooplus_k G_0(V_kG/V_{k-1}G)$
\end{lemme}

\begin{proof}
The last two assertions are standard for Rees modules associated with a good filtration. Let us prove first that $\cR_{V}G_0$ has finite type over $\CC[u,w]$. By lemma \ref{lemfacile} we know that $V_kG\cap G_0=0$ for $k\ll0$. It is then enough to prove that $\cR_{V}G_0$ is generated over $\CC[u,w]$ by $\oplus_{k\leq k_0}(V_kG\cap G_0)u^k$ for a suitable choice of $k_0$. But by the same lemma we know that there exists $k_0$ such that $V_kG\cap G_0=V_{k-1}G\cap G_0+V_kG\cap G_{-1}$ for all $k\geq k_0$. Hence for all such $k$, any element $m_ku^k\in (V_kG\cap G_0)u^k$ can be written
\begin{eqnarray*}
m_ku^k=u\cdot (m_{k-1}u_{}^{k-1})+w\cdot (n_{k-1}u_{}^{k-1})
\end{eqnarray*}
with $m_{k-1},n_{k-1}\in V_{k-1}G\cap G_0$. The result is obtained by an easy decreasing induction on $k$.

\medskip
Once we know that $\cR_{V}G_0$ has finite type, we can argue as follows to obtain the freeness. Because $\CC[u,u_{}^{-1},w]\otimes \cR_{V}G_0$ is isomorphic to $G_0[u,u_{}^{-1}]$, it is free of rank $\mu$ over $\CC[u,u_{}^{-1},w]$. Moreover this also implies that $\cR_{V}G_0$ is torsion free over $\CC[u,w]$, because its torsion must be supported on $\{u=0\}$ by this property, and this is not possible since $\cR_{V}G_0$ is contained in $G_0[u,u_{}^{-1}]$. Last, the fibers of $\cR_{V}G_0$ at each point of $\{u=0\}$ have dimension $\mu$: by homogeneity, it is enough to check this at $u=w=0$, where the fiber is $\oplus_kV_k(G_0/G_{-1})/V_{k-1}(G_0/G_{-1})$ and has dimension $\mu$ because of lemma \ref{lemfacile}. We conclude that $\cR_{V}G_0$ is a finite type projective $\CC[u,w]$-module, so is free of rank $\mu$ by Quillen-Souslin.
\end{proof} 

\subsubsection*{End of the proof of proposition \protect\ref{propbbases}}
Consider the Rees module $\cR_{V}G\defin\oplus_{k\in\ZZ}V_kG\cdot u^k$. It defines a meromorphic connection on $\PP^1\times \CC$ relative to the projection on $\CC$, with poles along $\{0\}\times \CC$ and $\{\infty\}\times \CC$: it is a free $\CC[u,w,w_{}^{-1}]$-module equipped with a compatible action of $\partial_w$ (here, $\infty=\{w=0\}$). Inside of it $\cR_{V}G_0$ is a relative lattice. The restriction $i^*\cR_{V}G$ to $\{u=0\}$ is equal to $\spe^{V}\!G\defin\oplus_k(V_kG/V_{k-1}G)$. This is a connection on $\PP^1$ with regular singularities at $0$ {\em and} $\infty$, and containing the lattice $\spe^{V}\!G_0=i^*\cR_{V}G_0$. We shall first construct a complementary lattice for $\spe^{V}\!G_0$ in $\spe^{V}\!G$ and we shall then extend it as a relative complementary lattice for $\cR_{V}G_0$ in $\cR_{V}G$. Restricting this lattice to $u=1$ will give a complementary lattice for $G_0$ in $G$. We shall then apply proposition \ref{propequiv} to obtain the good basis.

\medskip
First we define $H_{\beta}^{k}$ for every $\beta\in\QQ$ in the following way: let $\beta=\alpha-\ell$ with $\alpha\in{}[0,1[$ and $\ell\in\ZZ$; denote by $v$ the action of $\tau$ on $\spe^{V}\!G$. Put $H_{\beta}^{k}=v^\ell H_{\alpha}^{k-\ell}$. This space is contained in $v^\ell\gr_{\alpha}^{V}G=\gr_{\beta}^{V}G$ and is isomorphic to $H_{\alpha}^{k-\ell}$. We shall view $\gr_{\beta}^{V}G$ as a subspace of $V_{-\ell+1}G/V_{-\ell}G$ by the same arguments that we used for $E$ at the beginning of the proof.

For a fixed $\beta\in\QQ$, we have $H_{\beta}^{k}=0$ for $k\gg0$ and $H_{\beta}^{k}=\gr_{\beta}^{V}G$ for $k\ll0$: it is enough to prove this for $\beta=\alpha\in{}[0,1[$; this follows from the transversality condition and the fact that the analogue is true for the increasing filtration $G_k\gr_\alpha^{V}G$. It is easily seen that for every $\beta\in\QQ$, the decreasing filtration $H_{\beta}^{\bbullet}$ of $\gr_\beta^{V}G$ is opposite to $G_\bbullet\gr_{\beta}^{V}G$. Put
\begin{eqnarray*}
H^k&=&\ooplus_{\beta\in\QQ}H_{\beta}^{k}.
\end{eqnarray*}
One can see that $H^0$ is a $\CC[v]$-module of finite type. It is clearly stable under the action of $v\partial_v$ which is the action induced by $\tau\dtau$ on $\spe^{V}\!G$ because the $H_{\alpha}^{k}$ are so. It is free because it has no torsion, being contained in $\spe^{V}\!G$. Moreover, $H^k=v^k H^0$. We have
\begin{eqnarray}
H^0/vH^0&=&\ooplus_{\beta\in\QQ}H_{\beta}^{0}/H_{\beta}^{1}\nonumber\\
&=&\ooplus _{\beta\in\QQ}\gr_{H_\beta}^{0}(\gr_{\beta}^{V}G)\nonumber\\
&=&\ooplus_{\beta\in\QQ}\gr_{0}^{G}(\gr_{\beta}^{V}G) \quad\mbox{(transversality)}.\label{eqntransv}
\end{eqnarray}
This shows that $H^0$ has rank $\mu$ and that $\Ct\otimes_{\CC[\tau]}H^0=G$.

Hence $H^0$ and $\spe^{V}\!G_0$ can be glued together in a locally free sheaf on $\PP^1$ and because of the transversality assumption, they satisfy property (2) of proposition \ref{propequiv}. Consequently, this locally free sheaf is indeed free.

\medskip
Before going further on in the proof, let us remark that equality (\ref{eqntransv}) shows that the characteristic polynomial of the residue of the connection on $H^0$ is equal to the spectral polynomial of $(G,G_0)$. This residue is semisimple if and only if $\tau\dtau+\beta$ acts by $0$ on $H_{\beta}^{0}/H_{\beta}^{1}$ for all $\beta\in\QQ$, which is equivalent to the fact that $NH_{\alpha}^{k}\subset H_{\alpha}^{k+1}$ for all $\alpha\in{}[0,1[$ and all $k\in\ZZ$.

\medskip
Let us now come back to the proof. We shall extend the free sheaf that we have obtained above as a locally free sheaf on $\PP^1\times D$, where $D$ is a small neighbourhood of $u=0$ in $\CC$, in such a way that its restriction to $U_\infty\times D$ is equal to the restriction of $\cR_{V}G_0$ to this neighbourhood. This sheaf will be equipped with an action of $v\partial_v$ when restricted to $U_0\times D$, so is a relative logarithmic connection with pole along $\{0\}\times D$ on this open set. By restriction to $u=u_0\in D-\{0\}$ and using homogeneity in $u$, we shall obtain a locally free sheaf on $\PP^1$ extending $G_0$ and with a logarithmic connection with pole at $\{0\}$ on $U_0$. This sheaf is free, being a deformation of a free sheaf on $\PP^1$. We will also make the deformation in such a way that the residue at $0$ is constant in the deformation. So we shall obtain a good basis using proposition \ref{propequiv} and the previous remark.

\medskip
Using the isomorphism (\ref{eqnisom}), we can extend ${H^0}_{}^{\rm an}$ as a relative logarithmic connection $\cH^0$ in $\cR_{V}G_{|\Delta\times D}$, where $\Delta$ is a small neighbourhood of $v=0$ in $\PP^1$. Indeed, the isomorphism (\ref{eqnisom}) gives a trivialization $\cR_{V}G_{|\Delta\times D}\simeq \pi^*\spe^{V}\!G_{|\Delta}$, where $\pi:\PP^1\times D\rightarrow \PP^1$ denotes the projection. The glueing of ${H^0}_{}^{\rm an}$ and $\spe^{V}\!G_0$ on $\Delta-\{0\}$ can also be extended to a glueing of $\cH^0$ and $\cR_{V}G_0$ on $(\Delta-\{0\})\times D$, restricting $D$ if necessary: it is enough for this to extend bases of ${H^0}_{}^{\rm an}$ and $\spe^{V}\!G_0$ to the deformation. The last thing to verify in order that the program we have announced to be completed is that the residue of the relative connection on $\cH^0$ is constant: this is due to the fact that for all $\beta\in\QQ$ the graded deformation $\gr_\beta^{V}(\cH^0)$ is constant.
\end{proof}

\subsection*{Behaviour with respect to duality}
Assume that we are given an isomorphism
\begin{eqnarray*}
(G^*,G_0^*) \mbox{ or } (\ov{G^*},\ov{G_0^*})&\isom &(G,G_w).
\end{eqnarray*}
It defines a non degenerate bilinear (or sesquilinear) pairing
\begin{eqnarray*}
S:G_0\ootimes_{\CC[\theta]}^{}G_0&\longrightarrow &\theta^w\CC[\theta]
\end{eqnarray*}
of $\CC[\theta]$-modules which is compatible with the action of $\partial_\theta$. We also get a duality (see remark \ref{remdualspectre})
\begin{eqnarray*}
\gr_\alpha^VG\mbox{ or }\gr_\alpha^V\ov G&\isom&\Hom_\CC(\gr_{-\alpha+1}^{V}G,\CC)\quad (\alpha\in{}]0,1[)\\
\gr_0^VG\mbox{ or }\gr_0^V\ov G&\isom&\Hom_\CC(\gr_{0}^{V}G,\CC).
\end{eqnarray*}

A good (or very good) basis $\varepsilong$ of $G_0$ is said to be {\em adapted to} $S$ if
$$
S(\varepsilon_i,\varepsilon_j)\in\CC\cdot\theta^w\quad\mbox{for all }i,j.
$$

\begin{proposition}\label{propbbasesdual}
If the decreasing filtrations $H_\alpha^\bbullet$ of $\gr_\alpha^VG$ ($\alpha\in{}[0,1[$) satisfy conditions $1$ and $2$ of proposition {\rm\ref{propbbases}} and moreover
\begin{enumerate}
\refstepcounter{enumi}\refstepcounter{enumi}
\item
$H_{\alpha}^{k\perp}=H_{-\alpha+1}^{w-k}$ ($\alpha\in{}]0,1[$) and $H_{0}^{k\perp}=H_{0}^{w-k+1}$
\end{enumerate}
then the good basis given by this proposition is adapted to $S$.
\end{proposition}

\begin{proof}[Sketch of proof]
As in proposition \ref{propequiv}, one shows that the basis $\varepsilong$ is adapted to $S$ if and only if $S$ extends to a nondegnerate pairing
\begin{eqnarray*}
\cS:\cG_0\ootimes_{\cO_{\PP^1}^{}}^{}\cG_0&\longrightarrow &\cO_{\PP^1}^{}[w]
\end{eqnarray*}
compatible with connections, where $\cO_{\PP^1}^{}[w]$ denotes the trivial rank one bundle equipped with the connection $d+w\dpl\frac{d\theta}{\theta}$. Condition $3$ means that this is true for $\spe_VG_0$ and for the form $\spe_VS$ (see the proof of \ref{propbbases} for the notation). Now $\cR_VS$ defines $\cR_V\cS$ with maybe poles along $\{v=0\}\times D$. Since $\cR_V\cS_{|u=0}^{}=\spe_V\cS$ has no poles at $v=0$, it follows that $\cR_V\cS$ has no poles along $v=0$, and specializing to $u=1$ gives $\cS$.
\end{proof}

\section{Examples when there exists a good basis}
We first obtain from proposition \ref{propbbases}:
\begin{corollaire}
Let $(G,G_0)$ be such that no two distinct elements of the spectrum differ by an integer. Then $G_0$ has a good basis.
\end{corollaire}

\begin{proof}
In fact the assumption means that for every $\alpha\in{}[0,1[$ there exists a unique $k\in\ZZ$ (that we shall denote $k_\alpha$) such that $\gr_{k}^{G}\gr_{\alpha}^{V}G\neq0$, \ie. $G_{k_\alpha}\gr_{\alpha}^{V}G=\gr_{\alpha}^{V}G$ and $G_{k_\alpha-1}\gr_{\alpha}^{V}G=0$. So we can choose $H_{\alpha}^{k}=\gr_{\alpha}^{V}G$ for $k\leq k_\alpha$ and $H_{\alpha}^{k}=0$ for $k>k_\alpha$ in order to apply proposition \ref{propbbases}.
\end{proof}

\begin{rem*}
In fact, the $H_{\alpha}^{k}$ that we have constructed above is the only possible choice when the assumption of the corollary is satisfied. In particular the good basis that we obtain is very good if and only if $\tau\dtau+\alpha$ acts by $0$ on $\gr_\alpha^{V}G$ for every $\alpha\in{}[0,1[$, which means that the monodromy of $G$ is semisimple.
\end{rem*}

Let $c$ be a complex number and $G\otimes e_{}^{c\tau}$ be the $\CC[\tau,\tau^{-1}]$-module $G$ where the action of $\dtau$ is translated by $c$: $\dtau(g\otimes e_{}^{c\tau})=[(\dtau+c)g]\otimes e_{}^{c\tau}$. In the same way the action of $t=\theta^2\dth$ is translated by $-c$. We shall denote $G_0\otimes e_{}^{c\tau}$ the corresponding lattice.

\begin{proposition}
If $(G,G_0)$ admits a good basis $\bmg$ then $\bmg\otimes e_{}^{c\tau}$ is a good basis for $(G\otimes e_{}^{c\tau},G_0\otimes e_{}^{c\tau})$ and the matrix $A_1$ is the same for both\ptbl
\end{proposition}

\begin{proof}
The matrix of $t$ in the basis $\bmg\otimes e_{}^{c\tau}$ is equal to $A_0+c\id+\theta A_1$. In order to verify that $\bmg\otimes e_{}^{c\tau}$ is a good basis, it is enough to verify that the spectrum of $(G\otimes e_{}^{c\tau},G_0\otimes e_{}^{c\tau})$ is equal to the one of $(G,G_0)$ (the latter being equal to the spectrum of $A_1$ by assumption). It is then enough to verify that $V(G\otimes e_{}^{c\tau})=(VG)\otimes e_{}^{c\tau}$. This is clear because $\tau\dtau$ acts on $G\otimes e_{}^{c\tau}$ as $\tau\dtau+c\tau$ on $G$ and $c\tau$ acts by $0$ on each $\gr_\alpha^{V}G$.
\end{proof}

It follows from proposition \ref{proptensor} that

\begin{proposition}
If $(G',G'_0)$ and $(G'',G''_0)$ satisfy the assumptions of corollary {\rm\ref{corsymspectre}} and have good (or very good) bases, then the tensor product of these bases is a good (or very good) basis of $(G'\otimes G'',G'_0\otimes G''_0)$. \hfill\qed
\end{proposition}

\refstepcounter{equation}\label{rembbases}
\subsection{M. Saito's criterion}
M. Saito gives the following criterion for the existence of a very good basis (\cf. \cite[prop\ptbl 3.7]{MSaito89}): If for every $\alpha\in[0,1[$ and all $j>0$ the map $N^j:\gr_{\alpha}^{V}G\rightarrow \gr_{\alpha}^{V}G$ strictly shifts the filtration induced by $G_\bbullet$ by $j$, then all the conditions of proposition \ref{propbbases} are satisfied and hence $G_0$ admits a very good basis.

This criterion is fulfilled in particular if $\oplus _{\alpha\in[0,1[}^{}G_\bbullet\gr_{\alpha}^{V}G$ is the Hodge filtration of a mixed Hodge structure for which the weight filtration $M_\bbullet$ satisfies $NM_\bbullet\subset M_{\bbullet-2}^{}$.

Assume moreover that we have an isomorphism $(\ov{G^*},\ov{G_0^*})\isom(G,G_w)$. This isomorphism hence induces a perfect pairing of mixed Hodge structures
\begin{eqnarray*}
\ooplus_{\alpha\in{}]0,1[}^{} \gr_{\alpha}^{V}G\otimes \ooplus_{\alpha\in{}]0,1[}^{} \gr_{\alpha}^{V}G&\longrightarrow &\CC(-w)\\
\gr_0^VG\otimes \gr_0^VG&\longrightarrow &\CC(-w+1)
\end{eqnarray*}
where $\CC(k)$ is the pure Hodge of weight $-2k$ on $\CC$. Then M. Saito (see \cite[lemma 2.8]{MSaito89}) gives a canonical choice for the filtrations $H^\bbullet_\alpha$ (satisfying \ref{propbbases}-1,2), and this choice satisfies property 3 in proposition \ref{propbbasesdual}.

\section{Good basis for the Mellin transform}
\subsection*{Good basis and irregularity}
Let $G$ be as in section \ref{sec1}. Assume that $G$ admits a basis $\bmg$ over $\CC[\theta,\theta_{}^{-1}]$ for which the matrix of $t=\theta^2\partial_\theta$ is $A_0+\theta A_1$, where $A_0$ and $A_1$ are constant matrices. By assumption, $G$ is regular at $\tau=0$. Let $\ir_\infty(G)$ be the Malgrange-Komatsu irregularity number (see \eg. \cite{Malgrange74}) of $G$ at $\tau=\infty$ (\ie. $\theta=0$) and let $\mu=\rg G$. From the classical results of Turrittin and Katz on irregular singularities (see \eg. \cite{Malgrange91} or \cite{Varadarajan91}) on gets: 

\begin{proposition}\label{propinversible}
The matrix $A_0$ is invertible if and only if $\ir_\infty(G)=\mu$. \hfill\qed
\end{proposition}

\subsection*{Mellin transform of $G$}
Put $\sigma=-\theta\partial_\theta$ and identify the ring $\CC[\theta,\theta_{}^{-1}]\langle\partial_\theta\rangle$ with the ring of finite difference operators $\CC[\sigma]\langle\theta,\theta_{}^{-1}\rangle$, where we have the relation $\theta \sigma=(\sigma+1)\theta$. The (rational) Mellin transform $\gG$ of $G$ is the $\CC(\sigma)$-vector space $\CC(\sigma)\otimes_{\CC[\sigma]}G$, equipped with a structure of a left $\CC(\sigma)\langle\theta,\theta_{}^{-1}\rangle$-module. It is known to be a finite dimensional $\CC(\sigma)$-vector space. More precisely, for $G$ as in section \ref{sec1}, we have

\begin{proposition}\label{propdim}
$\dim_{\CC(\sigma)}\gG=\ir_\infty(G)$.
\end{proposition}

\begin{proof}
Let $G\otimes L_\alpha$ be the $\CC[\theta,\theta_{}^{-1}]$-module $G$ on which $\theta\partial_\theta$ acts by $\theta\partial_\theta+\alpha$, for $\alpha\in\CC$. Then $\dim_{\CC(\sigma)}\gG=\chi(\pDR (G\otimes L_\alpha))$ for $\alpha$ general enough, where the RHS is the Euler characteristics of the algebraic de~Rham complex of $G\otimes L_\alpha$ on $\CC^*$. The local index theorem \cite{Malgrange74} and the fact that $\chi(\pDRa (G\otimes L_\alpha))=0$ imply that $\chi(\pDR (G\otimes L_\alpha))=\ir_\infty(G\otimes L_\alpha)$, and the equality $\ir_\infty(G\otimes L_\alpha)=\ir_\infty(G)$ is easy.
\end{proof}

\subsection*{Good basis for the Mellin transform}
Let $G$ as above and assume that there exists a basis $\bmg$ of $G$ in which the matrix of $t$ is $A_0+\theta A_1$.

\begin{proposition}
We have $\dim_{\CC(\sigma)}\gG=\rk G$ if and only if $A_0$ is invertible. If this is satisfied, then $\bmg$ is also a $\CC(\sigma)$-basis of $\gG$, and the matrix of $\tau=\theta_{}^{-1}$ is $-A_{0}^{-1}(A_1+\sigma\id)$.
\end{proposition}

\begin{proof}
The first part is a consequence of propositions \ref{propinversible} and \ref{propdim}. We have $-\sigma\bmg=A_0\tau\bmg+A_1\bmg$. If $A_0$ is invertible, $\CC(\sigma)\cdot\bmg$ is a $\CC(\sigma)$-subspace of $\gG$ stable by $\tau$, the matrix of $\tau$ being $-A_{0}^{-1}(A_1+\sigma\id)$. It is also stable by $\tau_{}^{-1}$, and consequently $\bmg$ generates $\gG$ over $\CC(\sigma)$ because by assumption it generates it over $\CC(\sigma)\langle\tau,\tau_{}^{-1}\rangle$. Using the equality of dimensions we conclude that $\bmg$ is a $\CC(\sigma)$ basis of $\gG$.
\end{proof}

\begin{rem*}
If $\bmg$ is a {\em good} basis (definition \ref{defbb}) for $(G,G_0)$ and if $\dim_{\CC(\sigma)}\gG=\rg G$, we see that the determinant of $\tau$ in the $\CC(\sigma)$-basis $\bmg$ of $\gG$ satisfies
\begin{eqnarray*}
\det(\tau;\bmg)=\star \det(\sigma\id+A_1)=\star\SP_\psi(G,G_0;\sigma)
\end{eqnarray*}
where $\star$ is a nonzero constant.
\end{rem*}

\subsection*{Mellin transform and Fourier transform}
Assume that $G$ is the Fourier transform of a regular holonomic $\CC[t]\langle\partial_t\rangle$-module $M$ on which the action of $\partial_t$ is invertible (see \T\ref{sec2}). Let $M[t_{}^{-1}]$ be the localized module and $\gM$ its Mellin transform: we put $s=-t\partial_t$ and $\gM=\CC(s)\otimes_{\CC[s]}M$. This is a $\CC(s)\langle t,t_{}^{-1}\rangle$-module which has finite dimension over $\CC(s)$. Let $M_0\subset M$ be a finite $\CC[t]$-module stable by $\partial_{t}^{-1}$ and generating $M$ over $\CC[t]\langle\partial_t\rangle$. It defines a lattice $G_0$ of $G$ (see section \ref{sec2}). Let $\gG$ be the Mellin transform of $G$ as above. We can identify $\gG$ and $\gM$ using the isomorphism $\CC(s)\langle t,t_{}^{-1}\rangle\simeq\CC(\sigma)\langle \theta,\theta_{}^{-1}\rangle$ given by $s\mapsto\sigma+1$ and $t\mapsto-\theta\sigma=-(\sigma+1)\theta$. The previous results give the following:

\begin{proposition}\label{propmellin}\mbox{ }
\begin{enumerate}
\item
We have $\rg G=\sum_{c\in\CC}\mu_c(M)$ and $\dim_{\CC(s)}\gM=\sum_{c\in\CC^*}\mu_c(M)$. We have $\dim_{\CC(\sigma)}\gG=\rg G$ if and only if $0$ is not a singular point of $M$.
\item
Assume that there exists a good basis $\bmg$ for $(G,G_0)$. Then, if $0$ is not a singular point for $M$, $\bmg$ is a $\CC(s)$-basis for $\gM$ and the matrix of $t$ in this basis is $s(A_1+s\id)_{}^{-1}A_0$. Moreover the determinant of $t$ in this basis has poles (counted with multiplicities) on the opposite spectrum of $(G,G_0)$ exactly. \hfill\qed
\end{enumerate}
\end{proposition}

\part{The Gauss-Manin system \mbox{of a cohomologically tame function}}

\section{Some properties of cohomologically tame functions on an affine manifold}\label{secproptame}
Let $U$ be a smooth affine quasi-projective variety (over $\CC$) of dimension $n+1$ and $f:U\to\Afu$ be a regular function. We say that $f$ is {\em cohomologically tame} if there exists a compactification of $f$ given by a commutative diagram
$$
\xymatrix{
U\ar@{^{(}->}[r]^-{j}\ar[rd]^-{f}& X\ar[d]^-{F}\\
&\AA^1
}
$$
where $X$ is quasi-projective and $F$ is proper, such that the complex $\bR j_*\QQ_U$ has no vanishing cycle at infinity, namely, for all $c\in\Afu$, the vanishing cycle complex $\phi_{F-c}^{}\bR j_*\QQ_U$ is supported in at most a finite number of points, all at finite distance, \ie. contained in $U$. We shall denote $\Sigma$ the set of critical values of $f$ and $\mu=\mu(f)$ the sum of the Milnor numbers of $f$ at its critical points.

\subsection*{Examples}
If $f:\AA_{}^{n+1}\rightarrow \Afu$ is ``tame''  (see \protect\cite[lemma 4.3]{Broughton88}), \ie. if there exists $\varepsilon>0$ and a compact $K\subset U$ such that $\norme{\partial f}\geq \varepsilon$ out of $K$, then $f$ is cohomologically tame (with respect to the standard partial compactification of the graph $X\subset\PP^n\times \Afu$ of $f$).

In fact, Parusi\'nski has shown \cite{Parusinski95} that $f:\AA_{}^{n+1}\rightarrow \Afu$ is cohomologically tame with respect to the standard partial compactification of the graph $X\subset\PP^n\times \Afu$ of $f$ if and only if it satisfies a condition weaker than tameness, called the Malgrange condition, saying that when $f(x)$ remains bounded, there exists $\varepsilon>0$ such that $\norme{x}\norme{\partial f(x)}\geq \varepsilon$ for $\norme{x}$ sufficiently large.

In \cite{N-Z92} is introduced the notion of M-tameness for such a polynomial, which gives global Milnor fibrations in big balls in $\CC_{}^{n+1}$. It is not clear that any M-tame polynomial is cohomologically tame, but the tame ones, or the ones satisfying Malgrange condition are both cohomologically tame and M-tame.

Fix coordinates on $\AA_{}^{n+1}$. If $f:\AA_{}^{n+1}\rightarrow \Afu$ is convenient and nondegenerate with respect to its Newton polyhedron at infinity, then $f$ is tame \cite{Broughton88}, and this example is considered with some details in \T\ref{secnondeg}.

Consider now examples where $U$ is different from $\AA_{}^{n+1}$. The first one is the case where $U$ is a curve, so $f$ is a meromorphic function on the complete curve $\ov U$.

Let $U=(\CC^*)_{}^{n+1}$ and let $f$ be a Laurent polynomial. Assume that $f$ is nondegenerate with respect to its Newton polyhedron $\Delta_\infty(f)$ and is convenient \cite{Kouchnirenko76}. Let $F:X\rightarrow \Afu$ be the partial compactification of $f$ given by considering the closure of the graph of $f$ in the product of $\Afu$ with the toric compactification of $U$ defined by $\Delta_\infty(f)$. Then it follows from \cite[lemma (3.4)]{D-L91} that $f$ is cohomologically tame with respect to $X$.

\begin{rem*}
Consider a closed embedding of $F:X\rightarrow \Afu$ in $p:\cY\times \Afu\rightarrow \Afu$, where $\cY$ is smooth and $p$ is the second projection. Let $\kappa:X\hookrightarrow \cY\times \Afu$ be the inclusion. Then $f$ is cohomologically tame if and only if the constructible complex $\ccF=\bR(\kappa\circ j)_*\CC_U$ is noncharacteristic with respect to $p$ along $X-U$, in other words its characteristic variety or microsupport (see \eg. \cite{K-S90}) $\car\ccF\subset T^*(\cY\times \Afu)$ (\ie. the one of the corresponding bounded complex of regular holonomic $\cD$-modules on $\cY\times \Afu$) satisfies $\car\ccF\cap \lefpar \cY\times T^*\Afu\rigpar \subset\cY\times \Afu$ over the points of $\kappa(X-U)$ (see \eg. \cite[prop-def\ptbl 1.1]{Parusinski95}).

Remark also that the noncharacteristic property is satisfied by $\bR j_*\QQ_U$ (or $\ccF$ as above) if and only if it is satisfied by each of its perverse cohomology sheaves.
\end{rem*}

{\em In the remaining of the paper we will assume for simplicity that $n\geq 1$, \ie. $\dim U\geq 2$.}

\medskip
The following theorem (which is essentially well known) contains the main properties of the direct images sheaves of the constant sheaf $\QQ_U$ under the map $f$.

\begin{theoreme}\label{thmimdir}
Let $f:U\rightarrow \Afu$ be a cohomologically tame polynomial with $n+1=\dim U\geq 2$. Then
\begin{enumerate}
\item
the complex $\bR f_*\pQQ_U$ (\resp. $\bR f_!\pQQ_U$) has nonzero perverse cohomology in degrees $k\in[-n,0]$ (\resp. $k\in[0,n]$) at most; for $k<0$ (\resp. for $k>0$) $\pbR^{k}\!f_*\pQQ_U$ (\resp. $\pbR_{}^{k}\!f_!\pQQ_U$) is the perverse constant sheaf of rank $\dim H_{}^{k+n}(U,\QQ)$ (\resp. $\dim H_{c}^{k+n}(U,\QQ)$) on the affine line $\Afu$;
\item
the kernel and the cokernel (in the perverse sense) of the natural morphism $\pbR^0\!f_!\pQQ_U\rightarrow \pbR^0\!f_*\pQQ_U$ are constant perverse sheaves (with the same rank) on $\Afu$.

\item
The complex $\bR f_*\QQ_U$ has cohomology in degrees $k\in[0,n]$ at most; for $k<n$, $\bR^kf_*\QQ_U$ is the constant sheaf of rank $\dim H_{}^{k}(U,\QQ)$ on the affine line $\Afu$ and $\bR^n\!f_*\QQ_U[1]$ is perverse and equal to $\pbR^0\!f_*\pQQ_U$; 
\item
the complex $\bR f_!\QQ_U$ has cohomology in degrees $k\in[n,2n]$ at most and $\bR_{}^{n}\!f_!\QQ_U=\cH_{}^{-1}(\pbR^0\!f_!\pQQ_U)$;
\item
if $n\geq 2$, one has an exact sequence
\begin{eqnarray*}
0\rightarrow \cK\rightarrow \bR_{}^{n}\!f_!\QQ_U\rightarrow \bR_{}^{n}\!f_*\QQ_U \rightarrow \cC\rightarrow \bR_{}^{n+1}\!f_!\QQ_U\rightarrow 0
\end{eqnarray*}
where $\cK$ and $\cC$ are constant sheaves.
\end{enumerate}
\end{theoreme}

\begin{rems}\label{rempoids}\mbox{ }
\begin{enumerate}
\item
For $n=1$, the last two statements have to be slightly modified.
\item
As we shall see below (\cf. \T\ref{remdemohodgebr}), the image (in the perverse sense) of the morphism $\pbR^0\!f_!\pQQ_U\rightarrow \pbR^0\!f_*\pQQ_U$ is isomorphic to a direct sum $\ccT\oplus \sigma_*\cL[1]$, where $\sigma:\Afu-\Sigma\hookrightarrow \Afu$ denotes the open inclusion, $\cL$ is the local system made of the cohomology spaces $H^n(f_{}^{-1}(t),\QQ)$ and $\ccT$ is supported on $\Sigma$. In other words, this image satisfies the conclusion of the decomposition theorem.
\end{enumerate}
\end{rems}

Let us first give some consequences of the tameness of $f$ on the nearby and vanishing cycle sheaves of $F$. Let $c\in\Afu$ and $i_{F_{}^{-1}(c)}^{}:F_{}^{-1}(c)\hookrightarrow X$ be the closed inclusion. We will also denote $j:f_{}^{-1}(c)\hookrightarrow F_{}^{-1}(c)$ the open inclusion induced by $j:U\hookrightarrow X$. Let also $i:X-U\hookrightarrow X$ denote the closed inclusion.

\begin{proposition}\label{propcom}
If $f$ is cohomologically tame (with respect to $X$) then the following properties are satisfied:
\begin{enumerate}
\item
for each $c\in\Afu$, $\phi_{F-c}^{}(j_!\QQ_U)=\phi_{F-c}^{}(\bR j_*\QQ_U)$, where $j_!$ is the extension by $0$ and
\begin{eqnarray*}
i_{}^{-1}\phi_{F-c}^{}(j_!\QQ_U)=i_{}^{-1}\phi_{F-c}^{}(\bR j_*\QQ_U)=0;
\end{eqnarray*}
\item
for each $c\in\Afu$, $\psi_{F-c}^{}(\bR j_*\QQ_U)=\bR j_*(\psi_{f-c}^{}\QQ_U)$ and $\psi_{F-c}^{}(j_!\QQ_U)=j_!(\psi_{f-c}^{}\QQ_U)$;
\item
for each $c\in\Afu$, $i_{F_{}^{-1}(c)}^{-1}(\bR j_*\QQ_U)=\bR j_*\QQ_{f_{}^{-1}(c)}$.
\end{enumerate}
\end{proposition}

\begin{proof}
(1) It is clear that the two sheaves coincide on $f_{}^{-1}(c)$. By Verdier duality we have $j_!\pQQ_U=\bD(\bR j_*\pQQ_U)$ and $\phip_{F-c}^{}(j_!\pQQ_U)\simeq \bD (\phip_{F-c}^{}(\bR j_*\pQQ_U))$ (see \cite{Brylinski86}). The two complexes $\phi_{F-c}^{}(j_!\QQ_U)$ and $\phi_{F-c}^{}(\bR j_*\QQ_U)$, being dual up to a shift, have the same support, which is contained in $f_{}^{-1}(c)$ by assumption of tameness, so these two complexes coincide.

(2) It is equivalent to show the analogous equalities using the perverse functor $\psip$ and the perverse sheaf $\pQQ_U$. Then the first equality is Verdier dual to the second one. From (1) we deduce that $\phip_{F-c}^{}(j_!\pQQ_U)=j_!(\phip_{f-c}^{}\pQQ_U)$ so the second equality is equivalent to $i_{F_{}^{-1}(c)}^{-1}(j_!\QQ_U)=j_!\QQ_{f_{}^{-1}(c)}$ and it is enough to prove that $i_{}^{-1}i_{F_{}^{-1}(c)}^{-1}(j_!\QQ_U)=0$. But this is clear because the functors $i_{}^{-1}$ and $i_{F_{}^{-1}(c)}^{-1}$ commute.

(3) is a direct consequence of (1) and (2).
\end{proof}

\begin{corollaire}\label{corcom}
For all $c\in\Afu$,
\begin{enumerate}
\item
the complex $\phi_{f-c}^{}\QQ_U$ has cohomology in degree $n$ at most,
\item
on the sheaf $\QQ_U$, the functors $\bR f_*$ and $\bR f_!$ commute with the functors $\psi_{f-c}^{}$, $\phi_{f-c}^{}$, $i_{f_{}^{-1}(c)}^{-1}$ and $i_{f_{}^{-1}(c)}^{!}$.
\end{enumerate}
\end{corollaire}

\begin{proof}
(1) We know that $\phip_{f-c}^{}\pQQ_U$ is a perverse sheaf on$f_{}^{-1}(c)$. The assumption of tameness implies that its support consits of a finite number of points. Hence it has cohomology in degree $0$ at most.

(2) Let us show the result for $\bR f_*$ and $\psi_{f-c}^{}$ for instance. Let $t$ be the coordinate on $\Afu$. We have
\begin{eqnarray*}
\psi_{t-c}^{}(\bR f_*\QQ_U)&=&\psi_{t-c}^{}(\bR F_*\bR j_*\QQ_U)\\
&=&\bR F_*(\psi_{F-c}^{}(\bR j_*\QQ_U)) \mbox{ ($F$ is proper)}\\
&=&\bR F_*(\bR j_*(\psi_{f-c}^{}\QQ_U)) \mbox{ (proposition \ref{propcom})}\\
&=&\bR f_*(\psi_{f-c}^{}\QQ_U).
\end{eqnarray*}
\end{proof}

\begin{proof}[Proof of theorem \protect\ref{thmimdir}]\mbox{ }

\begin{lemme}\label{lemimdir}
If $f$ is cohomologically tame and $n\geq 1$, then
\begin{enumerate}
\item
the complex $\bR f_*\QQ_U$ has cohomology in degrees $m\in[0,n]$ at most;
\item
for $m\in[0,n-1]$, the sheaf $\bR^m\!f_*\QQ_U$ is a constant sheaf on $\Afu$ of rank $\dim H^m(U,\QQ)$.
\end{enumerate}
\end{lemme}

\begin{proof}
The first point does not depend on the tameness assumption: this is Artin theorem (see for instance \cite[Prop\ptbl 10.3.17]{K-S90}). It can be easily shown with the tameness assumption, because, due to proposition \ref{propcom}, restriction to fibers causes no problem, so it is equivalent to the fact that each fiber $f_{}^{-1}(c)$ has cohomology in degrees $m\in[0,n]$ at most, which can be proved using Morse theory (see \eg. \cite[Th\ptbl 10.3.8]{K-S90}).

\medskip
For the second point, remark that for each $c\in\Afu$ we have an exact sequence
\begin{eqnarray*}
\cdots\rightarrow \cH^m(i_{c}^{-1}\bR f_*\QQ_U)\longrightarrow \cH^m(\psi_{t-c}\bR f_*\QQ_U)\longrightarrow \cH^m(\phi_{t-c}\bR f_*\QQ_U)\rightarrow \cdots
\end{eqnarray*}
and we have $\cH^m(\psi_{t-c}\bR f_*\QQ_U)=\psi_{t-c}\bR^m\!f_*\QQ_U$ because $\bR f_*\QQ_U$ has constructible cohomology on $\Afu$. Moreover we clearly have $\cH^m(i_{c}^{-1}\bR f_*\QQ_U)=i_{c}^{-1}\bR^m\!f_*\QQ_U$. From corollary \ref{corcom} we deduce that $\cH^m(\phi_{t-c}\bR f_*\QQ_U)=0$ for $m\neq n$. Hence for $m<n$ we have $\phi_{t-c}\bR^m\!f_*\QQ_U=0$, which implies that $\bR^m\!f_*\QQ_U$ is a locally constant sheaf on $\Afu$, hence a constant sheaf.

\medskip
The Leray spectral sequence with $E_{2}^{mq}=H^q(\Afu,\bR^m\!f_*\QQ_U)$ degenerates at $E_2$, as follows from the previous results and this gives the statement concerning the rank of $\bR^m\!f_*\QQ_U$ for $m<n$.
\end{proof}

\begin{rem*}
At this point, one can compute the generic rank of $\bR^n\!f_*\QQ_U$ (or $\bR^n\!f_!\QQ_U$): one uses the fact that for a constructible complex $\cF$ on $\Afu$, the cohomology sheaves of which are local systems on $\Afu-\Sigma$, one has, for $t_0\in\Afu-\Sigma$,
\begin{eqnarray*}
\chi(\Afu,\cF)&=&\chi(\cF_{t_0})-\sum_{c\in\Sigma}\chi(\phi_{t-c}^{}\cF).
\end{eqnarray*}
When applied to $\cF=\bR f_*\QQ_U$ (or $\bR f_!\QQ_U$), this gives
\begin{eqnarray*}
\rg \bR^n\!f_*\QQ_U&=&\mu+h^n(U,\QQ)-h_{}^{n+1}(U,\QQ).
\end{eqnarray*}
\end{rem*}

\begin{lemme}\label{lempervers}
We have $\pbR^k\!f_*\pQQ_U=\bR^{n+k}\!f_*\QQ_U[1]$ for any $k\in\ZZ$.
\end{lemme}

\begin{proof}
For each $c\in\Afu$ the complex $\pQQ_{f_{}^{-1}(c)}^{}\defin \QQ_{f_{}^{-1}(c)}^{}[n]$ is perverse and we have an exact sequence of perverse sheaves on $f_{}^{-1}(c)$
\begin{eqnarray}\label{eqnperv}
0\longrightarrow \pQQ_{f_{}^{-1}(c)}^{}\MRE{} \psip_{f-c}^{}\pQQ_U\MRE{\can}\phip_{f-c}^{}\pQQ_U\longrightarrow 0.
\end{eqnarray}
We deduce a long exact sequence, by taking direct images and perverse cohomology, and using the fact that $\psip$ and $\phip$ commute with $\bR f_*$ (corollary \ref{corcom}) and perverse cohomology:
\begin{eqnarray*}
\cdots\longrightarrow \pbR^m\!f_*\pQQ_{f_{}^{-1}(c)}^{}\MRE{} \psip_{t-c}^{}\pbR^m\!f_*\pQQ_U\MRE{\can}\phip_{t-c}^{}\pbR^m\!f_*\pQQ_U \longrightarrow \cdots
\end{eqnarray*}
where $\pbR^m\!f_*\pQQ_{f_{}^{-1}(c)}^{}=H_{}^{m+n}(f_{}^{-1}(c),\QQ)$. Using the fact that $H_{}^{m+n}(f_{}^{-1}(c),\QQ)=0$ for $m>0$ and (tameness) $\phip_{t-c}^{}\pbR^m\!f_*\pQQ_U=0$ for $m\neq0$, we conclude that for $m\neq0$ $\pbR^m\!f_*\pQQ_U$ are (locally) constant sheaves on $\Afu$ and we obtain an exact sequence
\begin{eqnarray*}
0\longrightarrow H_{}^{n}(f_{}^{-1}(c),\QQ)\MRE{} \psip_{t-c}^{}\pbR^0\!f_*\pQQ_U\MRE{\can}\phip_{t-c}^{}\pbR^0\!f_*\pQQ_U \longrightarrow 0
\end{eqnarray*}
which gives surjectivity of $\can$ for $\pbR^0\!f_*\pQQ_U$.

Now, for a perverse sheaf $\cF$ on $\Afu$, the fact that $\can$ is onto is equivalent to the fact that $\cF$ is an ordinary sheaf (near $c$). Since the perverse cohomology sheaves of $\bR f_*\QQ_U$ are ordinary sheaves up to a shift by $1$, we conclude that these are also the ordinary cohomology sheaves up to a shift by $1$.
\end{proof}

The statements of the theorem for $\bR f_*\QQ_U$ are now proved. The perverse statements for $\bR f_!\QQ_U$ follow by Verdier duality and the nonperverse ones are straightforward consequences. The last point of the theorem follows then from the second one. Let us now prove it.

Consider the complex $i_{}^{-1}\bR j_*\QQ_U$, which is the cone of $j_!\QQ_U\rightarrow \bR j_*\QQ_U$. Proposition \ref{propcom} shows that we have $\phi_{F-c}^{}(i_{}^{-1}\bR j_*\QQ_U)=0$ for each $c\in\Afu$. This implies that the perverse cohomology groups of $\bR F_*(i_{}^{-1}\bR j_*\QQ_U)$ are (locally) constant sheaves up to a shift, hence the ordinary cohomology groups also. Consequently the (perverse) cohomology groups of the cone of $\bR f_!\QQ_U\rightarrow \bR f_*\QQ_U$ are constant sheaves, which implies the first part of (2).
\end{proof}

\section{The Gauss-Manin system of a regular function}\label{secGM}
Let $\cO_U$ be the sheaf of regular functions on the affine manifold $U$ and $f:U\rightarrow \Afu$ be a regular function. The Gauss-Manin system of $f$ is the complex $f_+\cO_U$ of $\cD_{\afu}$-modules.

Denote $\Omega^k(U)$ the space of differential forms of degree $k$ with polynomial coefficients. Then

(1) $H^k(\Omega^{^\bullet}(U),d)=H^k(U,\CC)$,

(2) if $f$ has only isolated critical points, the complex $(\Omega^{^\bullet}(U),df\wedge)$ has cohomology in degree $n+1$ only.

\medskip
The following is proved as in \cite{Malgrange75b} or \cite{Pham79} (the statement about regularity is well known):

\begin{proposition}\label{proppham}
$f_+\cO_U$ is represented by the complex $(\Omega^{\bbullet+n+1}(U)[\partial_t],d_f)$, where the differential $d_f$ is defined by $d_f(\sum_i\omega_i\partial_{t}^{i})=\sum_id\omega_i\partial_{t}^{i}-\sum_idf\wedge\omega_i\partial_{t}^{i+1}$. The cohomology modules $\cH^j(f_+\cO_U)$ are naturally equipped with a structure of a $\CC[t]\langle\partial_t\rangle$-module which makes them holonomic modules, regular even at infinity. Moreover $\cH^j(f_+\cO_U)=0$ for $j\not\in{}[-n,0]$. \hfill\qed
\end{proposition}

\begin{rems*}\mbox{ }
\begin{enumerate}
\item
We identify here algebraic $\cD_{\afu}$-modules with modules over $\CC[t]\langle\partial_t\rangle$. We always have $\cH_{}^{-n}(f_+\cO_U)=\CC[t]$.
\item
If $U=\AA_{}^{n+1}$, the (left) action of $\partial_t$ is invertible on $\cH^j(f_+\cO_U)$ for $-n<j\leq 0$. 
\end{enumerate}
\end{rems*}

The complex $(\Omega^{\bbullet+n+1}(U)[\partial_t],d_f)$ comes equipped with an increasing filtration $M_\bbullet$ defined by
\begin{eqnarray*}
M_k\Omega^{\bbullet+n+1}(U)[\partial_t]=\Omega^{\bbullet+n+1}(U)[\partial_t]_{\leq k+\bbullet}^{}
\end{eqnarray*}
where the filtration on the RHS is the one by the degree in $\partial_t$. This defines a filtration on the cohomology groups of this complex.

Assume from now on that $f$ is cohomologically tame.

\begin{proposition}\label{propminimal}
Under this assumption, for $j<0$, the $\CC[t]$-module $\cH^j(f_+\cO_U)$ is free of rank $\dim H_{}^{j+n}(U,\CC)$.
\end{proposition}

\begin{proof}
By the comparison theorem we have
$$\pDRa\cH^j(f_+\cO_U)=\pbR^j\!f_*\pDRa\cO_U=\pbR^j\!f_*\pCC_U$$
and from theorem \ref{thmimdir} we know that for $j<0$ this complex is the constant sheaf (up to a shift) of the right rank.
\end{proof}

\begin{corollaire}
If $U=\AA_{}^{n+1}$, the subcomplex $(df\wedge\Omega_{U}^{\bbullet},d)$ of $(\Omega_{}^{\bbullet+1},d)$ has cohomology in degree $0$ and $n$ only and the relative de~Rham complex $(\Omega_{U/\afu}^{\bbullet},d)$ has cohomology in degree $0$ and $n$ at most.
\end{corollaire}

\begin{proof}
Let $\eta\in\Omega^k(U)$ with $0< k<n$ and $\omega=df\wedge\eta$ such that $d\omega=0$, \ie. $df\wedge d\eta=0$. This implies that $d_f\omega=0$, hence $\omega=d_f\xi$ with $\xi\in\Omega^k(U)[\tau]$. Thanks to the de~Rham lemma, one may assume that $\xi\in\Omega^k(U)$ so that the previous equality is equivalent to $\omega=d\xi$ and $df\wedge\xi=0$. Put $\xi=df\wedge \eta'$. Then $\omega=d(df\wedge\eta')$.

The second part of the corollary is now clear.
\end{proof}

We shall denote $M=\cH_{}^{0}(f_+\cO_U)$. It is henceforth a holonomic $\CC[t]\langle\partial_t\rangle$-module with regular singularities (even at $t=\infty$). Let $\wh M$ be its Fourier transform (see section \ref{sec2}) and let $G=\wh M[\tau_{}^{-1}]=\wwwh{M[\partial_{t}^{-1}]}$. Then $G$ is a free $\CC[\tau,\tau_{}^{-1}]$-module of finite rank $\mu$.

\begin{rems}\label{remstrict}\mbox{ }
\begin{enumerate}
\item
Consider on $U\times \Afuc$ (where $\Afuc=\spec\CC[\tau]$) the $\cD_{U\times \afuc}^{}$-module $\cO_{U\times \afuc}^{}\cdot e_{}^{-\tau f}$. Then the complex $(\Omega^{\bbullet+n+1}(U)[\tau],d_f)$, where $d_f$ is defined as in proposition \ref{proppham}.(1), represents the direct image by $p:U\times \Afuc\rightarrow \Afuc$ of $\cO_{U\times \afuc}^{}\cdot e_{}^{-\tau f}$. In other words we have $\wh M=\cH^0(p_+\cO_{U\times \afuc}^{}\cdot e_{}^{-\tau f})$.
\item
The localized complex $(\Omega^{\bbullet+n+1}(U)[\tau,\tau_{}^{-1}],d_f)=p_+\cO_{U\times \afuc}^{}[\tau_{}^{-1}]\cdot e_{}^{-\tau f}$ has cohomology in degree $0$ only when $f$ has only isolated critical points, and this cohomology is equal to $G$.
\end{enumerate}
\end{rems}

\section{The Brieskorn lattice $G_0$ and its spectrum}

Let $M_0\subset M$ be the image of $\Omega_{}^{n+1}(U)$ in $M=\Omega_{}^{n+1}(U)[\partial _t]/d_f\Omega_{}^{n}(U)[\partial _t]$.

\begin{theoreme}\label{thmbriesk}
Assume that $f$ is cohomologically tame. Then $M_0$ is a free $\CC[t]$-module of finite rank generating $M$ over $\CC[t]\langle\partial_t\rangle$.
\end{theoreme}

From proposition \ref{propfourier} we conclude, putting $\theta=\tau_{}^{-1}$:

\begin{corollaire}
Let $G_0$ be the $\CC[\theta]$-module generated by $M_0$ in $G$. Then $G_0$ is a lattice in $G$ (see {\rm \T\ref{sec1}}) and the rank of $G$ is equal to the sum of the Milnor numbers of $f$ at its critical points. Moreover the fiber $G_0/\theta G_0$ at $\theta=0$ (\ie. $\tau=\infty$) is equal to $\Omega_{}^{n+1}(U)/df\wedge\Omega^n(U)$.
\end{corollaire}

\begin{rems}\label{remmonodromie}\mbox{ }
\begin{enumerate}
\item
If $U=\AA_{}^{n+1}$, we have $M=M[\partial _{t}^{-1}]$, $M_0$ is stable by $\partial _{t}^{-1}$ and, using the identification $G=\wh M=M$, $G_0$ is identified with $M_0$, viewed as a $\CC[\theta]$-module and not as a $\CC[t]$-module. This lattice $G_0$ is usually called the {\em Brieskorn lattice}. We have $G_0/\theta G_0\simeq\CC[x_0,\ldots ,x_n]/(\partial f/\partial x_0,\ldots ,\partial f/\partial x_n)$.
\item
If $U\neq\AA_{}^{n+1}$, then $G_0$ is the lattice obtained from $M_0$ after saturation by $\partial _{t}^{-1}$. It is the image of $\Omega_{}^{n+1}(U)[\tau_{}^{-1}]$ in $\Omega_{}^{n+1}(U)[\tau,\tau_{}^{-1}]/d_f\lefpar \Omega_{}^{n+1}(U)[\tau,\tau_{}^{-1}]\rigpar$. So not all the elements of $G_0$ are represented by differential forms: they are represented by polynomials in $\theta$ with coefficients in $\Omega_{}^{n+1}(U)$.
\item
If $\prod(X+\beta)_{}^{\nu_\beta}$ denotes the spectral polynomial of $(G,G_0)$ as defined in {\rm\T\ref{sec1}}, it follows from proposition \ref{proppolcar} and  \cite[cor\ptbl 1.13]{Bibi96a} that $\prod_\beta(T-\exp2i\pi\beta)_{}^{\nu_\beta}$ is the characteristic polynomial of the monodromy at infinity of $f$ on the cohomology $H^{n+1}(U,f_{}^{-1}(t);\QQ)$.
\end{enumerate}
\end{rems}

\begin{proof}[Proof of the {$\CC[t]$} finiteness of $M_0$]
Consider first the following situation: $\cY$ is a smooth projective variety, $\cU$ is a Zariski open set of $\cY$ and $\cZ=\cY-\cU$. Let $j:\cU\times \Afu\hookrightarrow \cY\times \Afu$ and $k:\cY\times \Afu\hookrightarrow \cY\times \PP^1$ denote the inclusions. Let $\cM$ be a regular holonomic $\cD_{\cU\times \afu}^{}$-module (regularity along $\cY\times \PP^1-\cU\times \Afu$ is included), and let $\cN\subset\cM$ be a coherent sub-$\cD_{\cU\times \afu/\afu}^{}$-module. The cohomology of $\bR j_*\cN$ admits a natural structure of a $\cD_{\cY\times \afu/\afu}^{}$-module. Here the differential operators are algebraic. The corresponding analytic objects are indicated with the exponent ``an''.

\begin{proposition}\label{proprelcoh}
Assume that the analytic de~Rham complex $\DR_{}^{\rm an}j_+\cM$ (which has constructible cohomology on $\cY\times \Afu$) has no vanishing cycle in some (analytic) neighbourhood of $\cZ\times \Afu$ with respect to the projection $\pi:\cY\times \Afu\rightarrow \Afu$. Then there exists a coherent $\cD_{\cU\times \afu/\afu}^{}$-submodule $\cN$ of $\cM$ satisfying $\cD_{\cU\times \afu}\cdot\cN=\cM$ and such that the cohomology of $\bR j_*\cN$ is coherent over $\cD_{\cY\times \afu/\afu}^{}$.
\end{proposition}

\begin{lemme}\label{lemdrholreg}
Let $\wt\cM$ be a regular holonomic $\cD_{\cY\times \PP^1}^{}$-module and let $\Sigma\subset\cY\times \Afu$ be the set of points $(y,c)$ such that $y\in\supp\phi_{\pi-c}^{}(\DR_{}^{\rm an}\wt\cM)$. There exists a coherent $\cD_{\cY\times \Pu/\Pu}^{}$-submodule $\wt\cN$ of $\wt\cM$ such that $\cD_{\cY\times \Pu}^{}\cdot\wt\cN=\wt\cM$ and 
\begin{eqnarray*}
\wt\cN[*(\cY\times \infty)]&=&\wt\cM[*(\cY\times \infty)]\quad\mbox{on } \cY\times \Pu-\Sigma.
\end{eqnarray*}
\end{lemme}

In particular $\wt\cM[*(\cY\times \infty)]$ is $\cD_{\cY\times \Pu/\Pu}^{}[*(\cY\times \infty)]$-coherent when restricted to $\cY\times \Pu-\Sigma$.

\begin{proof}
By GAGA it is enough to prove the result in the analytic category.

There exists only a finite number of critical values $c\in\PP^1$ for which $\phi_{\pi-c}^{}(\DR_{}^{\rm an}\wt\cM)$ is nonzero. Hence $\Sigma$ is contained in the union of a finite number of fibres of $\pi$ and is compact.

Let $c\in\PP^1$ and $V^{(c)}\wt\cM$ be the Malgrange-Kashiwara filtration of $\wt\cM$ along $\cY\times \{c\}$ (see \eg. \cite{Bibi87,M-S86,MSaito86}). It is known that each step $V^{(c)}_\alpha\wt\cM$ is relatively coherent in an analytic neighbourhood of $\cY\times \{c\}$ (this can be proved using resolution of singularities \cite[th\ptbl 4.12.1]{M-S86}, see also \cite[prop\ptbl 5.1.5]{Bibi87}). We put $\wt\cN=V^{(c)}_1\wt\cM$ in this neighbourhood.

It follows from the correspondence between regular holonomic $\cD$-modules and vanishing cycles (see \eg. \cite{M-S86}) that we have $\wt\cM=V^{(c)}_1\wt\cM=\wt\cN$ in a neighbourhood of $(y,c)\not\in\Sigma$.

We hence construct $\wt\cN$ by glueing the various $V^{(c)}_1\wt\cM$ in neighbourhoods of critical values with $\wt\cM$ outside these critical values.

For any critical value $c\in\PP^1$, let $z$ be a local coordinate on $\PP^1$ centered at $c$. We also have $\wt\cM[z_{}^{-1}]=(V^{(c)}_1\wt\cM)[z_{}^{-1}]$ (see \eg. \cite{M-S86}). Applying this for $c=\infty$ gives the result.
\end{proof}

\begin{proof}[Proof of proposition \protect\ref{proprelcoh}]
Let $\cM$ be an algebraic $\cD_{\cU\times \AA^1}^{}$-module which is holonomic and regular even at infinity. If $\DRa j_+\cM$ has no vanishing cycle in some neighbourhood of $\cZ\times \Afu$, then the same is true for $\DRa\cH^\ell(j_+\cM)$ for each $\ell$: in fact $\phip_{\pi-c}^{}\pDRa\cH^\ell(j_+\cM)$ is the $\ell$-th perverse cohomology object of $\phip_{\pi-c}^{}\pDRa(j_+\cM)$, and the last complex is zero if and only if all its perverse cohomology objects are so. Consequently the set $\Sigma$ associated with $j_+\cM$ is the union of the sets $\Sigma_\ell$ associated with $\cH^\ell(j_+\cM)$. By assumption we have $\cY\times \Pu=(\cU\times \Afu)\cup(\cY\times \Pu-\Sigma)$.

We will apply the lemma to $\wt\cM=\cH^0j_+\cM$. We get a coherent $\cD_{\cY\times \Pu/\Pu}^{}$-module $\wt\cN\subset\wt\cM$ and we put $\cN=j^*k^*\wt\cN$, which is a coherent $\cD_{\cU\times \afu/\afu}^{}$-module generating $\cM$ over $\cD_{\cU\times \afu}^{}$. Moreover on $\cY\times \Pu-\Sigma$ we have
\begin{eqnarray*}
(\bR k_*\bR j_*\cN)_{}^{\rm an}&=&(\bR (k\circ j)_*(k\circ j)^*\wt\cN)_{}^{\rm an}\\
&=&(\bR (k\circ j)_*(k\circ j)^*k_*k^*\wt\cN)_{}^{\rm an}\\
&=&(\bR (k\circ j)_*(k\circ j)^*k_*k^*\wt\cM)_{}^{\rm an}\mbox{ because of \ref{lemdrholreg}}\\
&=&(\bR k_*\bR j_*\cM)_{}^{\rm an}\\
&=&(k_+j_+\cM)_{}^{\rm an}.
\end{eqnarray*}
The last complex has $\cD_{\cY\times \Pu/\Pu}^{\rm an}[*(\cY\times \infty)]$-coherent cohomology on $\cY\times \Pu-\Sigma$ because of lemma \ref{lemdrholreg}. We conclude that $(\bR k_*\bR j_*\cN)_{}^{\rm an}$ is $\cD_{\cY\times \Pu/\Pu}^{\rm an}[*(\cY\times \infty)]$-coherent on $\cY\times \Pu$ and $\bR k_*\bR j_*\cN$ is $\cD_{\cY\times \Pu/\Pu}^{}[*(\cY\times \infty)]$-coherent by GAGA. Consequently $\bR j_*\cN$ is $\cD_{\cY\times \afu/\afu}^{}$-coherent.
\end{proof}

\begin{corollaire}\label{correlcoh}
Let $\cM$ be as in proposition {\rm \ref{proprelcoh}}. Then for any coherent $\cD_{\cU\times \afu/\afu}^{}$-submodule $\cN\subset\cM$ and for any $\ell$ the image of
\begin{eqnarray*}
\bH^\ell(\cU\times \Afu,\DR_\pi\cN)&\longrightarrow & \bH^\ell(\cU\times \Afu,\DR_\pi\cM)
\end{eqnarray*}
has finite type over $\CC[t]$, where $\DR_\pi$ denotes the algebraic de~Rham complex relative to $\pi:\cU\times \Afu\rightarrow \Afu$.
\end{corollaire}

\begin{proof}
It is enough to prove this for some $\cN$ generating $\cM$ over $\cD_{\cU\times \Afu}^{}$. One uses the $\cN$ given by the previous proposition. Then
$$
\bR(\pi\circ j)_*\DR_\pi\cN=\bR\pi_*\bR j_*\DR_\pi\cN= \bR\pi_*\DR_\pi\bR j_*\cN
$$
has $\cO_{\afu}$-coherent cohomology and
\begin{eqnarray*}
\bH^\ell(\cU\times \Afu,\DR_\pi\cN)&=&\Gamma(\Afu,\bR^\ell(\pi\circ j)_*\DR_\pi\cN).
\end{eqnarray*}
\end{proof}

\subsubsection*{End of proof of the finiteness of $M_0$}
Let $X\supset U$ be the quasi-projective variety associated with $f$ (see beginning of section \ref{secproptame}), let $\cX$ be a projective compactification of it and $\cY$ be a smooth projective manifold containing $\cX$ as a closed subset. Let $\cZ=\cX-U$ and $\cU=\cY-\cZ$. We have maps
\begin{eqnarray*}
\begin{array}{ccccccc}
U&\HMRE{i}&U\times \Afu&\HMRE{\eta}&\cU\times \Afu&\HMRE{j}&\cY\times \Afu\\
&&\MDL{p}&&\MDL{\pi}&&\MDL{\pi}\\
&&\Afu&=&\Afu&=&\Afu
\end{array}
\end{eqnarray*}
where $i$ and $\eta$ are closed immersions. We have $i_+\cO_U=\cO_U[\partial_t]\cdot \delta(t-f)$ and the relative de~Rham complex $\DR_\pi(i_+\cO_U)$ is the one given in proposition \ref{proppham}.

\medskip
It is enough to prove that if $N$ is a $\cD_{U\times \afu/\afu}^{}$-submodule of $i_+\cO_U$, the image of
\begin{eqnarray*}
H^\ell(\DR_p N(U\times \Afu))\longrightarrow H^\ell(\DR_p i_+\cO_U(U\times \Afu))
\end{eqnarray*}
has finite type over $\CC[t]$ for all $\ell$ where $\DR_\pi$ denotes the de~Rham complex relative to $p$: indeed, if we take for $N$ the $\cD_{U\times \afu/\afu}^{}$-submodule generated by $\cO_U\cdot \delta(t-f)$, we see that $M_0$ is contained in the image of
\begin{eqnarray*}
H_{}^{n+1}(\DR_p N(U\times \Afu))\longrightarrow H_{}^{n+1}(\DR_p i_+\cO_U(U\times \Afu)).
\end{eqnarray*}

The module $\eta_+i_+\cO_U$ satisfies the assumption of proposition \ref{proprelcoh} and since $\eta$ is a closed relative immersion, we can apply corollary \ref{correlcoh} to $\cN=\eta_+N$.
\end{proof}

\begin{proof}[Proof of the {$\CC[t]$} freeness of $M_0$]
Let $c$ be any critical value of $f$ and $V_\bbullet^{(c)} M$ be the Malgrange-Kashiwara filtration relative to $t-c$. It is enough to prove that, in a neighbourhood of $c$, we have $M_0\subset V_{<1}^{(c)}M$, because by construction $V_{<1}^{(c)}M$ has no torsion near $c$. Let $V_{\bbullet}^{(c)}(i_+\cO_U)$ be the corresponding Malgrange-Kashiwara filtration on $i_+\cO_U$.

\begin{lemme}
After restriction to the complement of all critical values $\neq c$, we have
\begin{eqnarray*}
V_\bbullet M&=&{\rm image}\;\cH^0\lefpar \DR_\pi(V_\bbullet(i_+\cO_U))\rigpar \rightarrow \cH^0\lefpar \DR_\pi(i_+\cO_U)\rigpar
\end{eqnarray*}
where $\pi:U\times \Afu\rightarrow \Afu$ denotes the projection.
\end{lemme}

\begin{proof}
As indicated in the proof of lemma \ref{lemdrholreg}, when we restrict to the complement $W_c$ of all critical values different from $c$, \ie. when we tensorize with $\CC[t,(\prod_{c'\neq c}^{}(t-c'))_{}^{-1}]$, each $V_\alpha (i_+\cO_U)$ is $\cD_{U\times W_c/W_c}^{}$-coherent. If $V'_\alpha M$ denotes the RHS in the lemma, we conclude that each $V'_\alpha M$ is finite over $\CC[t,(\prod_{c'\neq c}^{}(t-c'))_{}^{-1}]$. The other characteristic properties of the Malgrange-Kashiwara filtration are clearly satisfied by $V'_\bbullet M$, so $V'_\bbullet M=V_\bbullet M$ on $W_c$ by uniqueness.
\end{proof}

In order to conclude, it is enough to prove that $\delta(t-f)\in V_{<1}^{(c)}(i_+\cO_U)$. This is equivalent, by a standard argument, that the roots of the Bernstein for $(f-c)^s$ are negative. It is shown in \cite{Br-M90} (see also \cite[Prop\ptbl 4.2.1]{M-N91}) that the Bernstein polynomial for $(f-c)^s$ is equal to the lcm of the local analytic Bernstein polynomials of $f-c$ at the critical points of $f$ with critical value $c$. The roots of each of these local Bernstein polynomials are negative \cite{Kashiwara76}, so the same is true for the global Bernstein polynomial.
\end{proof}

\section{Duality for the Brieskorn lattice and for the spectrum} \label{secdualBriesk}
We will adapt the construction of higher residue pairings given by K. Saito \cite{KSaito83} in the present algebraic situation. We will construct in this section an isomorphism $\ov G^*\isom G$ which strictly shifts the filtration $G_\bbullet$ by $n+1$, \ie. $\ov G^*_{0}\isom G_{n+1}$. We then deduce from corollary \ref{corsymspectre} the following:

\begin{corollaire}\label{corsym}
The spectrum at infinity of a cohomologically tame polynomial of $n+1$ variables is symmetric with respect to $\dpl\frac{n+1}{2}$. \hfill\qed
\end{corollaire}

Proposition \ref{propresolstricte} shows that it is enough to construct a morphism
\begin{eqnarray*}
DM&\longrightarrow &M
\end{eqnarray*}
where $M=\cH^0f_+\cO_U$ and $D$ is the duality functor for $\CC[t]\langle\partial_t\rangle$-modules, such that the kernel and the cokernel are free $\CC[t]$-modules of finite type, and which strictly shifts by $n+1$ the microdifferential Brieskorn lattice $\muM_0$ of $f$.

Because $f_+\cO_U$ has cohomology in degrees $-n$ and $0$ only we see that $Df_+\cO_U$ has cohomology in degrees $0$ and $n$ only. Moreover we have $\cH^0Df_+\cO_U=DM$, and $\cH_{}^{-n}f_+\cO_U$ and $\cH^nDf_+\cO_U=D\cH_{}^{-n}f_+\cO_U$ are free $\CC[t]$-modules of rank one. Hence it is enough to construct a morphism
\begin{eqnarray}\label{eqndualsimple}
Df_+\cO_U&\longrightarrow &f_+\cO_U
\end{eqnarray}
such that the cohomology of its cone consists only of free $\CC[t]$-modules of finite type.

Let us consider a smooth quasi-projective compactification of $f$, namely a smooth quasi-projective manifold $\ov X$ and a commutative diagram
$$
\xymatrix{
U\ar@{^{(}->}[r]^-{j}\ar[rd]^-{f}&\ov X\ar[d]^-{\ov f}\\
&\AA^1
}
$$
with $\ov f$ proper and $j$ open. It will be convenient to assume that $D=\ov X-U$ is a divisor in $\ov X$ and that $\ov X$ dominates the partial compactification $X$ for which $f$ is cohomologically tame.

Because duality commutes with proper direct image (see \eg. \cite{Mebkhout87,MSaito89b,Schneiders86}), it is enough to construct a morphism $Dj_+\cO_U\rightarrow j_+\cO_U$. But we have such a canonical morphism because $j^+Dj_+\cO_U=\cO_U$. One verifies that such a morphism (\ref{eqndualsimple}) does not depend on the choice of the compactification.

By the comparison theorem and the local duality theorem (see \eg. \cite[th\ptbl 4.3.1]{Mebkhout87}) this morphism corresponds {\em via} the functor $\pDRa$ to the canonical morphism $\bR f_!\pCC_U\rightarrow \bR f_*\pCC_U$. But we know that the cohomology of the cone of this one has constant cohomology sheaves. Hence the cohomology of (\ref{eqndualsimple}) has $\CC[t]$-free cohomology modules. \hfill\qed

\subsection*{Relation with local duality}
Let $\cM=M_{}^{\rm an}=\cH^0\ov f_+(\cO_{\ov X}[*D])$. Let $B$ be an open neighbourhood of the critical points of $f$ in $U$ and $\Delta$ a neighbourhood of the critical values such that $f:B\rightarrow \Delta$ is the local Milnor fibration and let $\cM'=\cH^0f_+(\cO_{B}^{\rm an})$ be the local Gauss-Manin system of $f$ around its critical points. We have a natural restriction map $\cM_{|\Delta}^{}\rightarrow \cM'$ of $\cD_\Delta$-modules. Moreover we have a commutative diagram of Poincar\'e duality maps
$$
\begin{array}{ccc}
D\cM_{|\Delta}^{}&\MRE{}&\cM_{|\Delta}^{}\\
\MUR{}&&\MDR{}\\
D\cM'&\MRE{}&\cM'
\end{array}
$$
where the left vertical map is the adjoint of the right vertical one: this follows from the functoriality of the Poincar\'e duality map as constructed in \cite{Sch-Sch94} for instance.

It is known (see \eg. \cite{MSaito89}) that the local Poincar\'e duality map $D\cM'\rightarrow \cM'$ induces an isomorphism $D{\cM'}_{}^{\!\mu}\rightarrow {\cM'}_{}^{\!\mu}$ which sends $(D{\cM'}_{}^{\!\mu})_0$ onto ${\cM'}_{n+1}^{\!\mu}$ if ${\cM'}_{0}^{\!\mu}$ denotes the formal microlocal Brieskorn lattice of $f:B\rightarrow \Delta$ (in fact, this is true at the level of microlocal lattices). So it remains to identify ${\cM'}_{0}^{\!\mu}$ with $\muM_0$ in order to get the result, according to proposition \ref{propresolstricte}.

\subsection*{Identification of the microdifferential Brieskorn lattice}
Put $\cL=\cO_{\ov X}^{}[*D]$ and choose an $\cO_{\ov X}^{}$-coherent submodule $\cL_0$ of $\cL$ such that $\cD_{\ov X}^{}\cL_0=\cL$ and that the image of the composed map
$$
\ov{f}_*\lefpar \Omega_{\ov X}^{n+1}\ootimes_{\cO_{\ov X}^{}}^{}\cL_0\rigpar \longrightarrow \ov{f}_*\lefpar \Omega_{\ov X}^{n+1}\ootimes_{\cO_{\ov X}^{}}^{}\cL\rigpar=f_*\Omega_{U}^{n+1}\longrightarrow M
$$
is equal to the one of $f_*\Omega_{U}^{n+1}$, namely $M_0$. Such an $\cL_0$ exists since $M_0$ is $\cO_{\afu}^{}$-coherent.

Let $\cE_{\ov X}^{}$ be the sheaf of microdifferential operators on $T^*\ov X$ (formal or convergent, this will not matter now) with its subsheaf $\cE_{\ov X}^{}(0)$. Let $\cL_{0}^{\mu}=\cE_{\ov X}^{}(0)\otimes_{\cO_{\ov X}^{}}^{}\cL_0$.

\begin{lemme}
The image of the natural morphism
\begin{eqnarray*}
\bR^0\ov f_*\lefpar \cE_{\afua\leftarrow\ov X}^{}(0)\ootimes_{\cE_{\ov X}^{}(0)}^{\bL} \cL_{0}^{\mu}\rigpar &\longrightarrow & \bR^0\ov f_*\lefpar \cE_{\afua\leftarrow\ov X}^{}\ootimes_{\cE_{\ov X}^{}}^{\bL} \cL_{}^{\mu}\rigpar =\cE_{\afua}^{}\ootimes_{\cD_{\afua}^{}}^{}\cH^0\ov f_+\cL
\end{eqnarray*}
is equal to $\muM_0$.
\end{lemme}

\begin{proof}
This is a direct consequence from the fact that the natural map
\begin{eqnarray*}
\ov f_{}^{-1}\cE_{\afua}^{}(0)\ootimes_{\ov f_{}^{-1}\cO_{\afua}^{}}^{}\Omega_{\ov X}^{n+1}&\longrightarrow & \cE_{\afua\leftarrow\ov X}^{}(0)
\end{eqnarray*}
induces a quasi-isomorphism after $\otimes_{\cO_{\ov X}^{}}^{\bL}\cL_0$ and direct image by $\ov f$, as follows from \cite[\T4]{Kashiwara76}.
\end{proof}

Now $\cL_{}^{\mu}$ is supported on $D$ and on the critical points of $f$, and the part of the direct image coming from $D$ is zero, because $\ov f$ factorizes through $X$ where a non characteristic property is assumed along $X-U$. Hence $\muM_0$ is indeed the microlocal Brieskorn lattice of $f:B\rightarrow \Delta$. \hfill\qed

\section{The case of a convenient nondegenerate polynomial}\label{secnondeg}
Assume that $f:\AA_{}^{n+1}=U\rightarrow \Afu$ is nondegenerate with respect to its Newton polyhedron at infinity and that $f$ is convenient (see \cite{Kouchnirenko76}). Then it is known that $f$ is tame (see \cite{Broughton88}). One can define the Newton spectrum of $f$ (see \cite{Douai93}, but here we shall consider an increasing Newton filtration, so the Newton spectrum considered here is opposite to the one considered in \cite{Douai93}). We shall prove

\begin{theoreme}\label{thmnewtspec}
In this situation, the Newton spectrum of $f$ is equal to the spectrum of the Brieskorn lattice of $f$.
\end{theoreme}

\begin{rem*}
This result is analogous to the one of M. Saito \cite{MSaito88} (see also \cite{K-V85}) for the case of an isolated singularity. The proof will be analogous.
\end{rem*}

\begin{proof}
We shall use the following notation:

For a $n$-face $\sigma$ of the Newton polyhedron of $f$ not containing the origin, $L_\sigma$ will denote the linear form with coefficients in $\QQ_+^*$ such that $L_\sigma\equiv 1$ on $\sigma$.

For $u(x)\in\CC[x]$, we denote $\delta^*_\sigma(u)=\max_\nu L_\sigma(\nu+{\bf 1})$ where $\nu\in\NN^{n+1}$ is the exponent of a monomial in $u$. Moreover $\delta^*(u)$ will denote the maximum over all such $\sigma$ of $\delta^*_\sigma(u)$. We define $\delta_\sigma$ and $\delta$ in the same way, replacing $L_\sigma(\nu+{\bf 1})$ with $L_\sigma(\nu)$.

For $\alpha\in \QQ$ we put $\cN_\alpha\Omega_{}^{n+1}=\{u\cdot dx\mid \delta^*(u)\leq \alpha\}$. This defines an increasing filtration of $\Omega_{}^{n+1}(U)$ by finite dimensional vector spaces with $\cN_\alpha\Omega_{}^{n+1}=0$ for $\alpha\leq 0$.

We put $\cN_\alpha G_0={\rm image}\,\cN_\alpha\Omega_{}^{n+1}\subset G_0$.

We define a filtration $\cN_\alpha G$ of $G$:
\begin{eqnarray*}
\cN_\alpha G&=&\cN_\alpha G_0+\tau\cN_{\alpha+1}^{}G_0+\cdots+\tau^k\cN_{\alpha+k}^{}G_0+\cdots
\end{eqnarray*}
which is clearly stable under the action of $\CC[\tau]$ and satisfies $\tau\cN_\alpha G\subset \cN_{\alpha-1}^{}G$.

\begin{lemme}
The filtration $\cN_\bbullet G$ is equal to the Malgrange-Kashiwara filtration $V_\bbullet G$.
\end{lemme}

\begin{proof}
Let us first show that $\cN_\alpha G$ has finite type over $\CC[\tau]$ (over $\CC[\tau]\langle\tau\partial_\tau\rangle$ would be enough). It is enough to show that, for a given $\alpha$, there exists $k_0$ such that for $k\geq k_0$ we have
\begin{eqnarray*}
\tau^k\cN_{\alpha+k}^{}G_0\subset\CC[\tau]\cdot\lefpar \cN_\alpha G_0+\cdots+\tau_{}^{k_0}\cN_{\alpha+k_0}^{}G_0\rigpar .
\end{eqnarray*}
Let us fix $k_0$ such that $\cN_{\alpha+k_0}^{}G_0+G_{-1}^{}=G_0$, let $k\geq k_0$ and let $u \cdot dx\in\cN_{\alpha+k}^{}\Omega_{}^{n+1}(U)$. By the division lemma \cite[2.2.1]{Douai93}, we have
\begin{eqnarray*}
u\cdot dx&=&v\cdot dx+df\wedge\eta
\end{eqnarray*}
with $[v\cdot dx]\in\cN_{\alpha+k_0}^{}G_0$ and $\delta^*(d\eta)\leq \delta^*(u dx)-1$. So we have modulo $\im d_f$ the equality
\begin{eqnarray*}
\tau^k u\cdot dx&=&\tau^k v\cdot dx+\tau_{}^{k-1}d\eta
\end{eqnarray*}
and we can argue by decreasing induction on $k$ to get the result.

\medskip
Let $\sigma$ be a codimension one face not containing the origin and put $\xi_\sigma=L_\sigma(x\partial_x)$. Then $\lefpar \xi_\sigma+L_\sigma({\bf 1})\rigpar (ue_{}^{-\tau f})\cdot dx\in d\lefpar \Omega_{}^{n}(U)[\tau]e_{}^{-\tau f}\rigpar$ for $u\in\CC[x]$, so we have the following relation modulo $\im d$ in $\Omega_{}^{n+1}(U)[\tau]e_{}^{-\tau f}$:
\begin{eqnarray}\label{eqnnewt}
\qquad (\tau\partial_\tau+\delta^*_\sigma(u))\cdot ue_{}^{-\tau f}dx&\equiv&-[\xi_\sigma(u)-\delta_\sigma(u)u]e_{}^{-\tau f}dx+\tau(f-\xi_\sigma(f))ue_{}^{-\tau f}dx.
\end{eqnarray}
Moreover, for any such face $\sigma'$ we have $\delta^*_{\sigma'}\lefpar \xi_\sigma(u)-\delta_\sigma(u)u\rigpar \leq \delta^*_{\sigma'}(u)$ and $\delta_{\sigma'}^*((f-\xi_\sigma(f))u)\leq \delta_{\sigma'}^*(u)+1$, and both inequalities are strict if $\sigma'=\sigma$. We conclude that
\begin{eqnarray*}
\tau\partial_\tau\cdot \cN_\alpha G_0&\subset&\cN_\alpha G_0+\tau\cN_{\alpha+1}^{}G_0
\end{eqnarray*}
so $\cN_\alpha G$ is stable under the action of $\tau\partial_\tau$.

Moreover, by iterating $\#\{\sigma\mid \delta^*_\sigma(u)=\delta^*(u)\}$-times relation (\ref{eqnnewt}), one shows that there exists $N_0(\alpha)$ with $(\tau\partial_\tau+\alpha)_{}^{N_0(\alpha)}\cN_\alpha G_0\subset \cN_{<\alpha}^{}G$ and by the finiteness result above there exists $N(\alpha)$ such that $(\tau\partial_\tau+\alpha)_{}^{N(\alpha)}\cN_\alpha G\subset \cN_{<\alpha}^{}G$.

Because $\cN_\alpha G_0=0$ for $\alpha\leq 0$, we have $\cN_\alpha G=\tau\cN_{\alpha+1}^{}G$ for $\alpha\leq 0$.

Because the action of $\partial_\tau$ on $G_0$ is induced by the multiplication by $f$ on $\Omega_{}^{n+1}(U)$ and because $\delta(f)=1$, we have $\partial_\tau\cN_\alpha G_0\subset \cN_{\alpha+1}^{}G_0$, hence $\partial_\tau\cN_\alpha G\subset \cN_{\alpha+1}^{}G$. Moreover for $\alpha>0$ we have $\cN_{\alpha+1}^{}G=\cN_\alpha G+\partial_\tau\cN_{\alpha}^{}G$: indeed, if $ue_{}^{-\tau f}dx\in\cN_{\alpha+1}^{}G_0$, we have $(\partial _\tau\tau+\alpha)_{}^{N(\alpha+1)}ue_{}^{-\tau f}dx\in\cN_{<\alpha+1}^{}G$, so
\begin{eqnarray*}
ue_{}^{-\tau f}dx&=&\partial_\tau\lefcro \tau b(\tau\partial _\tau)ue_{}^{-\tau f}dx\rigcro +v
\end{eqnarray*}
with $v\in\cN_{<\alpha+1}^{}G$ and $\tau b(\tau\partial _\tau)ue_{}^{-\tau f}dx\in\cN_\alpha G$ and we iterate the process on $v$ to get the result.
\end{proof}

The spectrum of the Newton filtration is by definition the spectrum of the filtration $\cN_\bbullet$ defined by $$\cN_\bbullet(G_0/G_{-1}^{})= \cN_\alpha G_0/(\cN_\alpha G_0\cap G_{-1}^{}).$$ From the previous lemma we get $\cN_\alpha G_0\subset V_\alpha G\cap G_0$, hence $\cN_\alpha(G_0/G_{-1}^{})\subset V_\alpha(G_0/G_{-1}^{})$ for all $\alpha$. In order to show that both filtrations (or spectra) coincide, it is enough (by an argument of Varchenko, \cite{K-V85}) to show that both spectra are symmetric with respect to $(n+1)/2$. For $V$, this has been shown in the previous section and for $\cN$ this has been shown in \cite[prop\ptbl 7.3.3]{Douai93}, using arguments analogous to the ones in \cite{K-V85}.
\end{proof}

\begin{rems*}\mbox{ }
\begin{enumerate}
\item
R. Garcia pointed out that the relation between the Newton spectrum and the characteristic polynomial of the monodromy at infinity of $f$ that one deduces from the identification between the Newton spectrum and the spectrum of the Brieskorn lattice was expected in \cite{A-S89}.
\item
The order with respect to the Newton filtration of the $n+1$-form $dx$ is minimum and is obtained only for this form. It follows from \cite{Douai93} that the class of $dx$ in $G_0/G_{-1}^{}$ generates the vector space $V_{\alpha_{\min}^{}}^{}(G_0/G_{-1}^{})$, where $\alpha_{\min}^{}$ is the smallest spectral number. This space has thus dimension one.
\end{enumerate}
\end{rems*}

\section{Hodge theory for the Brieskorn lattice}\label{secappli}
We assume in this section that $f:U\rightarrow \Afu$ is a cohomologically tame regular function.

We will use the identification $$\psim_\tau\wh M=\psim_\tau(\cH^0\wwh{f_+\cO_U})=\psim_\tau(\cH^0\wwh{f_+\cO_U}[\tau_{}^{-1}])=\psim_\tau G=\ooplus_{\alpha\in[0,1[}^{}\gr_\alpha^VG$$
and it follows from from \cite[th\ptbl 5.3]{Bibi96a} that this space is equipped with a natural mixed Hodge structure isomorphic by \cite{MSaito87} to the limit of $H_{}^{n+1}(U,f_{}^{-1}(t))$ for $\tau\rightarrow \infty$ as constructed by Steenbrink and Zucker \cite{S-Z85} up to a Tate twist.

\begin{theoreme}\label{thmhodgebriesk}\mbox{ }
\begin{enumerate}
\item
The weight filtration of this mixed Hodge structure is the monodromy filtration of $N$ centered at $-n$ for $\psim_{\tau,\neq1}G$ and at $-(n+1)$ for $\psim_{\tau,1}G$.
\item
The Hodge filtration is the filtration induced by $G_\bbullet$ on $\psim_{\tau}G$.
\end{enumerate}
\end{theoreme}

The proof will be given in \T\ref{numproofthmhodgebriesk}. As a consequence of the proof we will obtain in \T\ref{numproofcorpositif} the positivity of the spectrum:

\begin{corollaire}\label{corpositif}
The spectral numbers of the Brieskorn lattice are contained in the closed interval $[0,n+1]$ and if $U=\AA_{}^{n+1}$ in the open interval $]0,n+1[$.
\end{corollaire}

\begin{rem*}
If $U\neq\AA_{}^{n+1}$ it may happen that $0$ belongs to the spectrum, as shown by the following example: take $U=(\CC^*)_{}^{n+1}$, $f$ is a Laurent polynomial which Newton polyhedron has dimension $n+1$ and contains $0$ in its interior (\ie. $f$ is convenient) and which is nondegenerate with respect to its Newton polyhedron. Then the class of the form $\dpl\frac{dx_0}{x_0}\wedge\cdots\wedge\frac{dx_n}{x_n}$ is contained in $G_0\cap V_0 G$.
\end{rem*}

From standard results in Hodge theory we have

\begin{corollaire}
The morphism $N:\psim_{\tau}G\rightarrow \psim_{\tau}G$ strictly shifts by one the filtration induced by $G_\bbullet$. \hfill\qed
\end{corollaire}

As a consequence we get from M. Saito's criterion \ref{rembbases}

\begin{corollaire}
The Brieskorn lattice of a cohomologically tame function has a very good basis $\varepsilong$ which satisfies $S(\varepsilon_i,\varepsilon_j)\in\CC\cdot\theta_{}^{n+1}$ if $S$ denotes the nondegenerate sesquilinear form $G\otimes_{\CC[\theta,\theta_{}^{-1}]}^{}G\rightarrow \CC[\theta,\theta_{}^{-1}]$ of {\rm\T\ref{secdualBriesk}}. \hfill\qed
\end{corollaire}

On deduces from proposition \ref{propequiv}

\begin{corollaire}
The Riemann-Hilbert-Birkhoff problem has a solution for the Bries\-korn lattice of a cohomologically tame function. \hfill\qed
\end{corollaire}

Let $V$ be any smooth quasi-projective manifold of pure dimension $\dim V$ and let $F_{D}^{\bbullet}H^k(V,\CC)$ be the (decreasing) Hodge-Deligne filtration on the cohomology spaces of $V$. For $q\in\NN$, let
\begin{eqnarray*}
\chi_{\rm Del}^{}(V;q)&\defin& (-1)_{}^{\dim V}\sum_k(-1)^k\dim\gr_{F_D}^{q}H^{k}(V,\CC)=\sum_i(-1)^i\dim\gr_{F_D}^{q}H^{\dim V-i}(V,\CC)\\
\zeta_{\rm Del}^{}(V;S)&\defin&\prod_q(S+\dim V-q)_{}^{\chi_{\rm Del}^{}(V;q)}.
\end{eqnarray*}
Remark that if we put $\chi_{\rm Del}^{c}(V;q)=\sum_i(-1)^i\dim\gr_{F_D}^{q}H_{c}^{\dim V+i}(V,\CC)$ we have $\chi_{\rm Del}^{c}(V;q)=\chi_{\rm Del}^{}(V;\dim V-q)$. For instance, if $V=\AA_{}^{n+1}$, we have $\zeta_{\rm Del}^{}(V;S)=(S+n+1)_{}^{(-1)_{}^{n+1}}$.

\begin{corollaire}\label{corspecU}
If $f$ is cohomologically tame on $U$, we have
\begin{eqnarray*}
\SP_\psi(G,G_0;S)&=&\SP_\phi(\psip_{1/t}^{}(\bR f_*\pCC_U);S)\cdot\zeta_{\rm Del}^{}(U;S).
\end{eqnarray*}
\end{corollaire}

\begin{proof}
As $f$ is cohomologically tame we know that $\cH^j(f_+\cO_U)$ are free $\CC[t]$-modules for $j\neq0$, so that $\cH^j(\wwh{f_+\cO_U})$ are supported at $\tau=0$ for $j\neq0$, and if we put $M=\cH^0(f_+\cO_U)$ we have
$$
\psim_\tau(\wwh{f_+\cO_U})=\psim_\tau(\cH^0\wwh{(f_+\cO_U)})=\psim_\tau\wh M=\psim_\tau G.
$$
Consequently we have $\SP_\psi(G,G_0;S)=\SP_\psi(\psip_\tau\wwwh{\bR f_*\pCC_U};S)$, according to the second part of theorem \ref{thmhodgebriesk}.

On the other hand, we get from \cite[cor\ptbl 5.4]{Bibi96a} that the RHS in \ref{corspecU} is equal to
$$
\SP_\phi(\phip_\tau\wwwh{\bR f_*\pCC_U};S)\cdot\zeta_{\rm Del}^{}(U;S)
$$
and we are reduced to showing
\begin{eqnarray*}
\SP_\psi(\psip_{\tau,1}\wwwh{\bR f_*\pCC_U};S)&=&\SP_\phi(\phip_{\tau,1}\wwwh{\bR f_*\pCC_U};S)\cdot\zeta_{\rm Del}^{}(U;S).
\end{eqnarray*}
This follows from the exact sequence in \cite[th\ptbl 4.3]{Bibi96a}. 
\end{proof}

{\em From now on, we will assume that $f:\AA_{}^{n+1}\rightarrow \Afu$ is a cohomologically tame polynomial.}

From the positivity statement \ref{corpositif} it follows that for a very good basis the matrix $A_1+k\id$ is invertible for all $k\in\NN$. From corollary \ref{corRH} we conclude

\begin{corollaire}
The partial Riemann-Hilbert problem for the Brieskorn lattice of $f$ has a solution. \hfill\qed
\end{corollaire}

Concerning the Aomoto complex considered in the introduction, the following corollary solves conjecture 7.4 in \cite{Douai93}. However, the problem of explicit computation of such a basis remains open in general, even in the case of a convenient nondegenerate polynomial.

\begin{corollaire}
If $0$ is not a critical value for $f:\AA_{}^{n+1}\rightarrow \Afu$, there exists a family of $\mu$ algebraic differential forms on $\AA_{}^{n+1}$ such that the determinant of the Aomoto complex computed in this basis is equal to $$c\cdot (s+1)^\mu/\SP_\psi(G,G_0;s+1)$$ with $c\in\CC^*$.
\end{corollaire}

The constant $c$ is equal to the product $\prod_if(x_{}^{(i)})_{}^{\mu_i}$ where $x_{}^{(i)}$ are the critical points of $f$ and $\mu_i$ the corresponding Milnor numbers (see \cite{L-S90}).

\begin{proof}
We may replace the differential $d_f$ of the complex $(\Omega_{}^{\bbullet+n+1}[1/f][\partial_t],d_f)$ with $d'_f\omega=d\omega-\dpl\frac{df}{f}\wedge\omega\cdot\partial_tt$, hence this complex isomorphic to $(\Omega_{}^{\bbullet+n+1}[1/f][s],d_s)$ with $s=-\partial_tt$ and $d_s=f_{}^{-s}\cdot d\cdot f^s$. This defines an isomorphism of the Mellin transform of $M[t_{}^{-1}]$ with $H^0(\Omega_{}^{\bbullet+n+1}[1/f](s),d_s)$ over the automorphism
\begin{eqnarray*}
\CC(s)\langle t,t_{}^{-1}\rangle&\longrightarrow &\CC(s)\langle t,t_{}^{-1}\rangle\\
\varphi(s,t)&\longmapsto&\varphi(s+1,t).
\end{eqnarray*}
The result is then a consequence of proposition \ref{propmellin}.
\end{proof}

We can restate cor\ptbl \ref{corspecU}:
\begin{corollaire}\label{corspec}
If $U=\AA_{}^{n+1}$ and $f$ is a cohomologically tame polynomial, we have $$\SP_\phi(\psi_{1/t}^{}\bR^n f_*\CC_U;S)=\SP_\psi(G,G_0;S).$$
\end{corollaire}

\begin{proof}
We have
\begin{eqnarray*}
\SP_\phi(\psip_{1/t}^{}\bR f_*\CC_U;S)&=& \SP_\phi(\psi_{1/t}^{}\bR^n f_*\CC_U;S)\cdot\SP_\phi(\psi_{1/t}^{}f_*\CC_U;S)_{}^{(-1)_{}^{-n}}
\end{eqnarray*}
and the computation made in \cite[rem\ptbl 5.5]{Bibi96a} shows that $\SP_\phi(\psi_{1/t}^{}f_*\CC_U;S)=S+n+1$. We may then apply corollary \ref{corspecU}.
\end{proof}

Consider now the operator of multiplication by $f$ on $G_0/\theta G_0=\Omega_{}^{n+1}/df\wedge\Omega^n$. It sends $V_\beta (G_0/\theta G_0)$ in $V_{\beta+1}^{}(G_0/\theta G_0)$ for each $\beta\in\QQ$, hence defines a nilpotent endomorphism $$[f]:\sum_\beta\gr_{\beta}^{V}(G_0/\theta G_0)\longrightarrow \sum_\beta\gr_{\beta+1}^{V}(G_0/\theta G_0)$$ (remark however that the multiplication by $f$ is not nilpotent on $G_0/\theta G_0$ in general). Let $T_\infty$ be the monodromy of $f$ along a positively oriented circle of big radius.

\begin{corollaire}
The nilpotent part of $T_{\infty}^{-1}$ and $[f]$ have the same Jordan structure.
\end{corollaire}

\begin{proof}
One may replace $T_{\infty}^{-1}$ with $\Tc_0$ (see \cite[th\ptbl 1.10]{Bibi96a}). The result is then proved as in \cite{Varchenko81} (see also \cite[\T7]{S-S85}).
\end{proof}

\begin{rem*}
We do not get here a precise relation between the nilpotent part of the monodromy at infinity and the operator of multiplication by $f$, as in the case of an isolated singularity.
\end{rem*}

\refstepcounter{equation}\label{numproofthmhodgebriesk}
\subsection{Proof of theorem \protect\ref{thmhodgebriesk}}
We take notation of \cite[\S\T4,5]{Bibi96a}: let $\cX$ be a smooth compactification of $U$ such that $\cX-U$ is a divisor with normal crossings and $f:U\rightarrow \Afu$ extends as $F:\cX\rightarrow \PP^1$. According to \cite[th\ptbl 5.3]{Bibi96a}, the first point would follow from the fact that
\begin{eqnarray*}
W_k\psip_\tau\cH^0\wwwh{f_+\cO_U}&=&\lefacc
\begin{array}{ll}
0&\mbox{if }k<0\\
\psip_\tau\cH^0\wwwh{f_+\cO_U}&\mbox{if }k\geq 0.
\end{array}
\rig
\end{eqnarray*}
We know that $W_k\psip_\tau\cH^0\wwwh{f_+\cO_U}=\psip_\tau\lefpar \wwwh{W_k\cH^0f_+\cO_U}\rigpar$ where $W_k\cH^0f_+\cO_U$ is the image of
$$
\cH^0F_+\lefpar W_kj_+\cO_U\rigpar \longrightarrow \cH^0F_+\lefpar j_+\cO_U\rigpar
$$
where $j:U\hookrightarrow F_{}^{-1}(\Afu)$ denotes the inclusion (see \cite[\S\T2.4.4 and 5]{Bibi96a}). It is then enough to show that the Fourier transform of $\cH^0(f_+\cO_U)/W_0\cH^0(f_+\cO_U)$ is supported at $\tau=0$.

It follows from \cite{MSaito87} that $(j_+\cO_U,W[\dim U]_\bbullet)$ underlies a mixed Hodge module as well as $(\cH^0(f_+\cO_U),W[\dim U]_\bbullet)$ (if the Hodge filtration of $\cO_U$ is such that $\gr_k^F\cO_U=0$ for $k\neq-\dim U$). Hence $\cH^0(f_+\cO_U)$ has weights $\geq \dim U$.

As a mixed Hodge module we have $\DD\cO_U=\cO_U(\dim U)$, hence $j_!\cO_U=\lefpar \DD j_+\cO_U\rigpar (-\dim U)$ so $j_!\cO_U$ has weights $\leq \dim U$ as well as $\cH^0f_!\cO_U$.

Consequently the morphism $\cH^0(f_!\cO_U)\rightarrow \cH^0(f_+\cO_U)$ factorizes through $W[\dim U]_{\dim U}^{}=W_0\cH^0(f_+\cO_U)$ and as the cokernel of this map is isomorphic to a power of $\CC[t]$ (tameness of $f$), the same is true for $\cH^0(f_+\cO_U)/W_0\cH^0(f_+\cO_U)$. \hfill\qed

\begin{rem}\label{remdemohodgebr}
This argument can be used to show remark \ref{rempoids}-(2): the image of $\cH^0(f_!\cO_U)\rightarrow \cH^0(f_+\cO_U)$ is a pure Hodge module.
\end{rem}

For the second point, recall that $\wh M=\cH^0\pch_+(\kappa_+\Etf)$ (see \loccit.). Let $G_\bbullet\kappa_+\Etf$ be the filtration defined and used in \cite[\T4]{Bibi96a}. It defines a filtration $G_\bbullet\pch_+(\kappa_+\Etf)$ of the complex $\pch_+(\kappa_+\Etf)$, hence a filtration
$$
G'_\bbullet \wh M=G'_\bbullet\cH^0\pch_+(\kappa_+\Etf)= \image\lefcro \cH^0\lefpar G_\bbullet\pch_+(\kappa_+\Etf)\rigpar \rightarrow \cH^0\pch_+(\kappa_+\Etf)\rigcro .
$$
It also defines a filtration $G_\bbullet\psim_\tau\kappa_+\Etf$ using the $V$-filtration (see \loccit.), and as above it defines a filtration
$$
G''_\bbullet\cH^0\pch_+(\psim_\tau\kappa_+\Etf)= \image \lefcro \cH^0\lefpar G_\bbullet\pch_+(\psim_\tau\kappa_+\Etf)\rigpar \rightarrow \cH^0\pch_+(\psim_\tau\kappa_+\Etf)\rigcro .
$$
As we have $\cH^0\pch_+(\psim_\tau\kappa_+\Etf)=\psim_\tau\cH^0\pch_+(\kappa_+\Etf)=\psim_\tau\wh M=\psim_\tau G$, we get two filtrations on $\psim_\tau G$: the filtration $G''_\bbullet\psim_\tau G$ is the one used in \cite[th\ptbl 5.3]{Bibi96a} and the filtration $G'_\bbullet\psim_\tau G$ which is naturally induced from $G'_\bbullet \wh M$.

Last, we denote $G_\bbullet$ the filtration obtained from the Brieskorn lattice $G_0$ and $G_\bbullet\psim_\tau G$ the filtration it induces naturally, related to the spectrum of the Brieskorn lattice.

We want to show that $G''_\bbullet\psim_\tau G=G_\bbullet\psim_\tau G$, in order to apply \cite[th\ptbl 5.3]{Bibi96a}.

\begin{lemme}\label{lemgp}
We have an inclusion $G'_\bbullet\psim_\tau G\subset G_\bbullet\psim_\tau G$.
\end{lemme}

\begin{proof}
We can filter the complex $(\Omega_{\cX}^{\bbullet+n+1}[*(D\cup F_{}^{-1}(\infty))]\otimes_\CC\CC[\tau],d-\tau df)$ either by the filtration counting the total pole order plus the degree in $\tau$, and after taking the global sections on $\cX$ one obtains $G'_\bbullet\wh M$, or by the filtration by the degree in $\tau$ only and one obtains the filtration of $\wh M$ by the $M_k=\sum_{j=0}^{k}\tau^jM_0$. Clearly the inclusion is true at the level of complexes, so it remains true at the cohomology level, hence $G'_\bbullet \wh M\subset M_\bbullet$ and $G'_\bbullet \psim_\tau\wh M\subset M_\bbullet\psim_\tau\wh M$.

On the other hand consider $G=\wh M[\tau_{}^{-1}]$. Then $G_k$ is defined as $\sum_{\ell\geq 0}^{}\tau_{}^{-\ell}\lefpar \mbox{image}M_k\rigpar$, so the isomorphism $\psim_\tau\wh M\rightarrow \psim_\tau G$ sends $M_\bbullet\psim_\tau\wh M$ in $G_\bbullet\psim_\tau G$.
\end{proof}

\begin{lemme}\label{lemgs}
We have $G''_\bbullet\psim_\tau G\subset G'_\bbullet\psim_\tau G$.
\end{lemme}

\begin{proof}
We will use the following
\begin{proposition}[{\cite[th\ptbl 4.8.1]{M-S86}}]
Let $\cM$ be a holonomic $\cD_\cX[\tau]\langle\partial_\tau\rangle$-module and let $V_\bbullet\cM$ be the Malgrange-Kashiwara filtration of $\cM$ relative to $\tau=0$. Then for all $i$ and all $\alpha$ we have a commutative diagram:
\begin{eqnarray*}
\begin{array}{ccc}
\cH^i\pch_+V_\alpha\cM&\MRE{}&V_\alpha\cH^i\pch_+\cM\\
\MDL{}&&\MDL{}\\
\cH^i\pch_+\gr_\alpha^V\cM&=&\gr_\alpha^V\cH^i\pch_+\cM
\end{array}
\end{eqnarray*}
where all the maps are onto. \hfill\qed
\end{proposition}

Then on the one hand we have
\begin{eqnarray*}
G''_\bbullet\gr_\alpha^VG&\defin&{\rm image}\;\cH^0G_\bbullet\pch_+\gr_\alpha^V(\kappa_+\Etf)\rightarrow \cH^0\pch_+\gr_\alpha^V(\kappa_+\Etf)
\end{eqnarray*}
and we will show below that
\begin{eqnarray}\label{eqnonto}
\cH^0G_\bbullet \pch_+V_\alpha(\kappa_+\Etf)\rightarrow \cH^0G_\bbullet\pch_+\gr_\alpha^V(\kappa_+\Etf)\quad\mbox{is onto}
\end{eqnarray}
so
\begin{eqnarray*}
G''_\bbullet\gr_\alpha^VG&=&{\rm image}\;\cH^0G_\bbullet\pch_+V_\alpha(\kappa_+\Etf)\rightarrow \cH^0\pch_+\gr_\alpha^V(\kappa_+\Etf).
\end{eqnarray*}

But on the other hand the image of $\cH^0G_\bbullet\pch_+V_\alpha(\kappa_+\Etf)$ in $\cH^0\pch_+(\kappa_+\Etf)$ is contained in $G'_\bbullet\cH^0\pch_+(\kappa_+\Etf)\cap V_\alpha\cH^0\pch_+(\kappa_+\Etf)$, hence we eventually get $G''_\bbullet\gr_\alpha^VG\subset G'_\bbullet\gr_\alpha^VG$ for all $\alpha\in[0,1[$.

\medskip
Let us show (\ref{eqnonto}). The computation of $(\gr_{\alpha}^{V}\kappa_+\Etf,G_\bbullet)$ made in \cite[\T4]{Bibi96a} shows that for $\alpha<1$, the isomorphism $\tau:\gr_{\alpha}^{V}\kappa_+\Etf\rightarrow \gr_{\alpha-1}^{V}\kappa_+\Etf$ sends $G_\bbullet$ onto $G_{\bbullet+1}^{}$. From the Hodge property the complex $G_\bbullet\pch_+\gr_\alpha^V(\kappa_+\Etf)$ is strict for any $\alpha\in{}[0,1]$, hence, using this isomorphism, it is strict for any $\alpha<1$. In particular we get that for any $\alpha\leq 1$,
\begin{eqnarray}\label{eqngrV}
\cH^kG_\bbullet\pch_+\gr_\alpha^V(\kappa_+\Etf)=0 \quad\mbox{for }k>0.
\end{eqnarray}

On the other hand, (\ref{eqnonto}) would follow from the fact that $\cH^1G_\bbullet\pch_+V_\alpha(\kappa_+\Etf)=0$ for any $\alpha<1$, and using (\ref{eqngrV}) it is enough to verify this for $\alpha\ll0$. So let us fix $p\in\ZZ$ and consider $\cH^1G_p\pch_+V_\alpha(\kappa_+\Etf)$. For a fixed $k$, we have $G_k(\kappa_+\Etf)\cap V_\alpha(\kappa_+\Etf)=0$ for $\alpha$ sufficiently small (see the definition of these filtrations in \loccit.), so for a fixed $p$ we have $G_p\pch_+V_\alpha(\kappa_+\Etf)=0$ for $\alpha\ll0$. \hfill\qed

\medskip
We hence have $G_\bbullet\gr_\alpha^VG\supset G''_\bbullet\gr_\alpha^VG$ and it is enough to prove that equality holds in order to get the second part of the theorem. This is done as in \cite[\T6]{S-S85} using the symmetry of the spectrum of $(G,G_0)$ on the one hand and the Hodge symmetry and the fact that the weight filtration is the one computed above on the other hand.
\end{proof}

\refstepcounter{equation}\label{numproofcorpositif}
\subsection{Proof of corollary \protect\ref{corpositif}} 
By definition, the filtration $G_\bbullet\kappa_+\Etf$ satisfies $G_k=0$ for $k<0$. It follows that the filtration $G''_\bbullet\psim_\tau(\kappa_+\Etf)$ considered above satisfies the same property and that the filtration that it induces on $\pDRa\psim_\tau(\kappa_+\Etf)$ by the formula (5.2) of \cite{Bibi96a} also. Hence $G''_k\psim_\tau G=0$ for $k<0$, so by the previous theorem we have $G_k\psim_\tau G=0$ for $k<0$. This implies in particular that $\nu_{\alpha+k}^{}=0$ for $\alpha\in[0,1[$ and $k<0$, so $\nu_\beta=0$ for $\beta< 0$. Using the symmetry of the spectrum \ref{corsym} one obtains the result.

If $U=\AA_{}^{n+1}$, then $G=\wh M$ and $G_0=M_0$. Lemmas \ref{lemgp} and \ref{lemgs} also hold for $\phim_\tau\wh M$ and we can conclude that $\nu_{\alpha+k}^{}=0$ for $\alpha\in[0,1]$ and $k<0$, so $\nu_\beta=0$ for $\beta\leq 0$. \hfill\qed

\providecommand{\bysame}{\leavevmode\hbox to3em{\hrulefill}\thinspace}

\end{document}